\documentclass[reqno]{amsart}
\usepackage{latexsym,amssymb,amsthm,amsmath}

\usepackage[bookmarksnumbered, colorlinks]{hyperref}
\usepackage{amssymb}
\usepackage{amsmath}
\usepackage{multicol}
\theoremstyle{plain}
\newtheorem{theorem}{Theorem}[section]
\newtheorem*{Theorem B}{Theorem B}
\newtheorem*{Theorem A}{Theorem A}
\newtheorem*{Theorem C}{Theorem C}
\newtheorem{lemma}[theorem]{Lemma}

\newtheorem{corollary}[theorem]{Corollary}

\newtheorem{definition}[theorem]{Definition}

\numberwithin{equation}{section}

\theoremstyle{remark}
\newtheorem{remark}[theorem]{Remark}
\begin{document}
\title[Li-Yau-Hamilton and Hessian estimates and applications]{Li-Yau, Hamilton  gradient and Hessian estimates for nonlinear weighted parabolic equations and applications}

\author[S.K. Hui, A. Abolarinwa, S. Bhattacharyya]{Shyamal Kumar Hui, Abimbola Abolarinwa and Sujit Bhattacharyya}


\begin{abstract}
	 This article is devoted to the study of several estimations for a positive solution to a nonlinear weighted parabolic equation on a   weighted Riemannian manifold.  We therefore derive new Li-Yau type and Hamilton type gradient estimates yielding several consequences. We also derive Hessian estimate and some corollaries for the same equation. Among the applications of our estimates discussed here are Harnack type inequalities,  Liouville type theorems and a local time reversed Harnack inequality. 
\end{abstract}

\maketitle

\section{Introduction}
Gradient and Hessian estimations are  effective tools to study nature of the solutions of linear and nonlinear partial differential equations(PDE).  In this method, bounds for the gradient and Hessian of a positive solution to a PDE can be found without solving the equation. This article is therefore devoted to finding some types of gradient and Hessian estimations for positive smooth solution of the following nonlinear parabolic equation
\begin{equation}\label{eq_heat_old}
	(\mathcal{L}_f-q-\partial_t)w = G(w),
\end{equation}
on  an $n$-dimensional weighted Riemannian manifold denoted by the triplet $(M^n,g,e^{-f}d\mu)$ endowed with metric $g$ and weighted volume form $e^{-f}d\mu$.  Here the functions  $q=q(x,t)$ and $f=f(x,t)$ ($x\in M$, $t\in (0,T]$) are space and time dependent smooth functions, that is $q,f \in C^\infty(M\times (0,T])$, whereas $0<w=w(x,t)$ denotes the space and time dependent solution to (\ref{eq_heat_old}).   The sufficiently smooth function $G :\mathbb{R}^+\to \mathbb{R}$ depending on $w>0$ depicts the nonlinear part of the model (\ref{eq_heat_old}) and includes several forms of nonlinearities.  To illustrate the generality and physical and geometric sources of (\ref{eq_heat_old}), we give different examples of $G$ soon. The compactness and completeness of the $M$ are suitably chosen.  As notations $\partial_t$ denotes the partial derivative with respect to $t$ (we can also write $\partial_t u =u_t$) and $\mathcal{L}_f$ denotes the weighted Laplacian (also called Witten Laplacian) defined by
$$\mathcal{L}_f u:= \Delta u-\langle \nabla f,\nabla u\rangle,$$ for any smooth function $u$ on $M$, where $\nabla$, $\Delta$  and $\langle \cdot,\cdot\rangle$ denote the usual gradient operator, Laplace-Beltrami operator and $g$ dependent inner product.
Putting $w=e^u$  ($u=\log w$) into (\ref{eq_heat_old}) we have
\begin{equation}\label{eq_heat_new}
	\mathcal{L}_f u =u_t-|\nabla u|^2+q+ \tilde{G}(u),
\end{equation}
where $\tilde{G}(u)=e^{-u}G(e^u)$, which is another nonlinear weighted parabolic equation.  Transforming (\ref{eq_heat_old}) to (\ref{eq_heat_new}) is a basic idea of Li-Yau gradient and Harnack estimates, where  (\ref{eq_heat_new}) motivates a definition of an appropriate Harnack quantity involving time derivative and norm of gradient of $u$ (where one obtains $u_t=w_t/w$ and $|\nabla u| =|\nabla w|^2/w^2$).  The Harnack quantity can be proved to evolve by a certain nonlinear heat flow which will enable the application of the maximum principle.  Integrating this Harnack quantity along smooth path in the underlying manifold yields some classical Harnack estimates.

The idea of Li-Yau and Hamilton type estimation was originated after the work of Li and Yau \cite{Li-Yau} for the equation
\begin{eqnarray}
	\nonumber (\Delta - q(x,t)-\partial_t)w(x,t) = 0,
\end{eqnarray} on complete Riemannian manifold and Hamilton \cite{Hamilton,Hamilton-harnack} for
\begin{eqnarray}
	\nonumber (\Delta -\partial_t)w(x,t) = 0,
\end{eqnarray} on compact Riemannian manifold. Note that setting $G(w)=0$ and $G(w)=q(x,t)=0$ respectively in (\ref{eq_heat_old}) give the last two equations, respectively. Thus, estimations on  (\ref{eq_heat_old})  should extend Li-Yau and Hamilton gradient estimates in \cite{Li-Yau} and  \cite{Hamilton,Hamilton-harnack}.  Later Souplet and Zhang \cite{SOUPLET-ZHANG} improved the estimations of Li and Yau,  yielding numerous interesting results including Liouville type theorems on non-compact Riemannian manifolds. Perelman's \cite{Perelman-1,Perelman-2,Perelman-3} revolutionary insights on the analytic and geometric aspect of Ricci flow popularized this method and a new perspective was given to study nonlinear partial differential equations on Riemannian manifolds.
This particular  field of research has a vast literature.
Yau's work then further extended by Li \cite{Li} for weighted Riemannian manifold involving weighted Laplacian ($f$-Laplacian). 

It is known that there is a direct connection between the $f$-Laplacian and the  Schr\"odinger operators in quantum mechanics and quantum field theory. If we consider $w$ to be a static solution i.e.,  $w$ does not evolve in time meaning that $\partial_t w=0$ then (\ref{eq_heat_old}) reduces to a generalized elliptic Poisson type equation. With certain choice of the function $G(w)$ we can get numerous applications for (\ref{eq_heat_old}) in Geometry, Physics, Chemistry and Biology. To illustrate the idea here is a short literature that shows the importance of studying (\ref{eq_heat_old}).
\begin{enumerate}
	\item[(a)] If we take $G(w)=aw^p-bw^q,\ a,b>0,\ q>p\ge 1$ then (\ref{eq_heat_old}) has a great importance in mathematical physics \cite{ABIMBOLA-7,DUNG-KHANH,Mastrolia}. A simple calculation would give 
	\begin{eqnarray}
		\nonumber \tilde{G}(w) &=& aw^{p-1} - bw^{q-1},\\
		\nonumber \nabla \tilde{G} &=& (a(p-1)w^{p-2}-b(q-1)w^{q-2})|\nabla w|,\\
		\nonumber \mathcal{L}_f\tilde{G} &=& \left(a(p-1)(p-2)w^{p-3}-b(q-1)(q-2)w^{q-3}\right)\mathcal{L}_fw.
	\end{eqnarray}
	\item [(b)] Caffarelli and Lin's \cite{Caffarelli} work on non-local heat flow preserving the $L^2$-norm involves the equation (\ref{eq_heat_old}) for $G(w)=\lambda(t)w+A(x,t)$ with the normalization condition $\int_M w^2 d\mu =1$, where $\lambda(t)=-\int_M (w\Delta w+wA(x,t))d\mu$. Again just like before we have  
	\begin{eqnarray}
		\nonumber \tilde{G}(w) &=& \lambda(t)+\frac{A(x,t)}{w},\\
		\nonumber \nabla \tilde{G} &=& \frac{\nabla A}{w} - \frac{A}{w^2}\nabla w,\\
		\nonumber \mathcal{L}_f\tilde{G} &=& \frac{\mathcal{L}_f(A(x,t))}{w}+A(x,t)\mathcal{L}_f\left(\frac{1}{w}\right) + \langle \nabla A(x,t), \nabla (\frac{1}{w})\rangle.
	\end{eqnarray}
	\item[(c)] The case when $G(w)=|w|^{b-1}w$ with $b>1$ has influence in Physics and Geometry, one can see \cite{Castorina,Castorina-2} for more information. In this case we have
	\begin{eqnarray}
		\nonumber \tilde{G}(w) &=& |w|^{b-1},\\
		\nonumber \nabla \tilde{G} &=& (b-1)|w|^{b-2}\nabla w,\\
		\nonumber \mathcal{L}_f\tilde{G} &=& (b-1)(b-2)|w|^{b-3}|\nabla w|^2 + (b-1)|w|^{b-2}\mathcal{L}_fw.
	\end{eqnarray}
	\item[(d)] Another important form of (\ref{eq_heat_old}) when $G(w)=aw(\log w)^\alpha$, $a,\alpha\in\mathbb{R}$, can be found in \cite{Polyanin}. For $\alpha=1$ the static part of the equation can be related with gradient Ricci solitons and Sobolev type inequalities with logarithmic nonlinearity \cite{DUNG-KHANH-2,Gheru}. Again we infer
	\begin{eqnarray}
		\nonumber \tilde{G}(w) &=& a(\log w)^\alpha,\ \nabla \tilde{G}\ =\ a\alpha\frac{(\log w)^{\alpha -1}}{w}\nabla w,\\
		\nonumber \mathcal{L}_f\tilde{G} &=& a\alpha (\alpha-1) \frac{(\log w)^{\alpha -2}}{w^2}|\nabla w|^2-a\alpha\frac{(\log w)^{\alpha -1}}{w^2}|\nabla w|^2\\
		\nonumber&&+a\alpha\frac{(\log w)^{\alpha -1}}{w}\mathcal{L}_f w. 
	\end{eqnarray}
	\item[(e)] Allen-Cahn type equation for $G(w)=cw(1-w^2),\ c>0,\ 0\le w<1$ is well known in the method of phase separation
	in iron alloys. This equation is well analyzed, one can see \cite{Calatroni} and the references therein for more information. Like previous we get
	\begin{eqnarray}
		\nonumber \tilde{G}(w) &=& c(1-w^2),\\
		\nonumber \nabla \tilde{G} &=& -2cw\nabla w,\\
		\nonumber \mathcal{L}_f\tilde{G} &=& -2c|\nabla w|^2 -2cw\mathcal{L}_fw.
	\end{eqnarray}
	\item[(f)] Last but not the least the well known Fisher-KPP equation \cite{Fischer,Fitz,Kolomogorov} for $G(w)=cw(1-w)$, $c>0,\ 0\le w<1$, has been studied in numerous fields of science, for example population density of genes \cite{Alza-Mano,Levin}, Fitzhugh-Nagumo equation for impulse growth in nerve axons \cite{Fitz}, flame propagation, nuclear reactors, combustion theory and heat transfer \cite{Luford,Tremel} etc. We also have
	Like previous we get
	\begin{eqnarray}
		\nonumber \tilde{G}(w) &=& c(1-w),\\
		\nonumber \nabla \tilde{G} &=& -2c\nabla w,\\
		\nonumber \mathcal{L}_f\tilde{G} &=& -2c\mathcal{L}_f w.
	\end{eqnarray}
\end{enumerate}
In the main result section, we have remarked (see Remark \ref{remark_important}) on the significance of the gradient and weighted Laplacian values for each case.\\

 In modern day, Brighton \cite{BRIGHTON} studied the Liouville type theorem on smooth metric measures space. Fang \cite{Fang} derived differential Harnack inequalities for heat type equations under geometric flow. In 2014, Zhu and Li \cite{Zhu-Li} established Li-Yau type estimates for nonlinear parabolic equations on manifolds. Elliptic gradient estimates for a weighted heat type equation  were studied by Wu \cite{WU-4}. He also studied gradient estimates for nonlinear parabolic equations \cite{WU-2, WU-3}. In the time span 2019-20, Abolarinwa et al. \cite{ABIMBOLA-4, ABIMBOLA-4.2,ABIMBOLA-JGP22} studied elliptic and parabolic gradient estimates on smooth metric measure space. Shen and Ding \cite{SHEN-DING} studied the blow-up conditions in porous medium systems. Dai et al. \cite{DAI} studied gradient estimates for weighted $p$-Laplacian equations on Riemannian manifolds. Zhao and Shen \cite{ZHAO-SHEN} studied gradient estimates for $p$-Laplacian on non compact metric measures space. Very recently, Hui et al. \cite{Sujit-4} established Li--Yau type estimation on evolving manifolds. One can see \cite{Sujit-1,Sujit-2,Sujit-3} and the references therein for further studying.

On the other hand Hessian estimation was being conducted simultaneously with gradient estimation. The Hessian of a $C^2$ function $f$ on a Riemannian manifold $(M,g)$ is a $(0,2)$ type tensor $\text{Hess }f:\chi(M)\times\chi(M)\to C^\infty(M)$ defined by $$\text{Hess }f(X,Y)=g(\nabla_X(\nabla f),Y),$$ where $\nabla_X$ denotes the covariant derivative with respect to $X$ and $\nabla f$ is the gradient of $f$. Taking trace of the above equation gives the Laplace-Beltrami operator (Laplacian operator) $$\Delta f=\text{div }(\nabla f).$$
The idea of Hessian estimation is similar to gradient estimation except for the case: we are going to find bound for the quantity $\frac{|\nabla ^2 w|}{w}$, where $\nabla^2$ denotes the Hessian of $w$ (a positive solution of (\ref{eq_heat_old})). This method was originated and popularized after the work of Hamilton \cite{Hamilton-harnack} in 1993, where he showed the matrix $A_{ij}=(\log u)_{ij}$, which is the Hessian matrix of $\log u$, satisfies $$A_{ij}\ge -\frac{1}{2t}g_{ij},$$ along the Ricci flow at time $t$. Later Han and Zhang \cite{Han-Zhang} studied upper bounds for Hessian matrices for positive solution of heat type equations on Riemannian manifolds. Li \cite{Yi-Li} extended the results of Han-Zhang to $V$-Laplacian, which is defined as $$\Delta_V u=\Delta u-\langle V,\nabla u\rangle, \text{ for some }V\in\chi(M).$$ As a special case of this operator if we set $V=\nabla f$, where $f$ is a smooth function on $M$, then $V$-Laplacian reduces to weighted Laplacian. Sharp bounds for Hessian of positive solution of Allen-Cahn type equation on Riemannian manifolds has been studied by Zhang \cite{Zhang}. In recent days, Abolarinwa et al.\cite{ABIMBOLA-6} studied Harnack inequalities for a class of nonlinear heat equation of the form $$\frac{\partial}{\partial t}w(x,t) = \Delta w(x,t) + g(x,t).$$ Local Hessian estimates of positive solutions to the non-linear parabolic equation $$\frac{\partial u}{\partial t} = \Delta u + \lambda u^\alpha$$
 have been studied by Wang et al. \cite{Wang-hessian}. This topic has been the subject of several research studies.
 
Motivated by such interesting applications and rich literature in this paper we derive Li-Yau type gradient estimation and Hamilton type estimation for positive solution of (\ref{eq_heat_old}) on a static manifold (i.e. the metric is independent time) and Hessian estimate for the same with some applications as corollaries. It should be mentioned here that we are going to use an important result during the estimation which is known as Shi's \cite{Shi} local derivative estimate. This will help us to bound the Riemann curvature and its derivative on the manifold. 
\begin{remark}
In any method of gradient and Hessian estimation it is quite natural to assume a bound for the weight function $f$ and its gradient $\nabla f$ but in this article we present estimations independent of such restriction.
\end{remark}
Section 2 consists of some basic results and informations. The section 3 is devoted to the statement and proofs of the results: first,  local and global Li-Yau type gradient estimate followed by local and global Hamilton type gradient estimate, and lastly, Hessian estimates. Their application to Harnack inequalities and Liouville type results are presented accordingly.  A brief summary of the results is given in the last section.
\section{Preliminaries}
\begin{definition}[\cite{Bakry-Emery}Bakry-\'Emery Ricci tensor]
	For any two positive integers $m\ge n$ and any smooth function $f$ on $M$, the $(m-n)$-Bakry-\'Emery Ricci tensor is given by
	\begin{eqnarray}
	\nonumber Ric_f^{m-n} := Ric+\text{Hess }f-\frac{\nabla f\otimes \nabla f}{m-n}.
	\end{eqnarray}
\end{definition}
\begin{remark}
	The case $m=n$ happens if and only if $f$ is a constant function. Furthermore, $m\to\infty$ gives the $\infty$-Bakry-\'Emery Ricci tensor 
	\begin{eqnarray}
	\nonumber Ric_f^\infty = Ric+\text{Hess }f,
	\end{eqnarray}
which is sometimes written as $Ric_f$.
\end{remark}
\begin{lemma}[Weighted Bochner Formula]\label{lemma_Bochner}
	For any smooth function $u$ on a weighted Riemannian manifold we have 
	\begin{eqnarray}\label{BF}
 \frac12 \mathcal{L}_f|\nabla u|^2 = |\text{Hess }u|^2+\langle \nabla \mathcal{L}_f u,\nabla u \rangle+Ric_f(\nabla u,\nabla u),
	\end{eqnarray}
where the norm on the Hessian is the usual Hilbert-Schmidt norm on $2$-tensors. 
\end{lemma}
Upon using Cauchy-Schwarz inequality and referring to the definition of weighted Laplacian $\mathcal{L}_f$,  we get
$$|\text{Hess }u|^2 \ge \frac{1}{n}(\Delta u)^2 = \frac{1}{n}(\mathcal{L}_f u+\langle \nabla f,\nabla u \rangle)^2.$$ 
Applying an elementary inequality  of the form 
$$(a+b)^2\ge \frac{a^2}{1+\alpha}-\frac{b^2}{\alpha},$$ where $\alpha>0$ with the choice $a=\mathcal{L}_f u,\ b=\langle \nabla f,\nabla u \rangle$ and $\alpha=\frac{m-n}{n}>0$ into (\ref{BF}),  we can get (\ref{WBI}) of the next lemma.
\begin{lemma}[Weighted Bochner inequality]
	For any smooth function $u$ on a weighted Riemannian manifold, we have
	\begin{eqnarray}\label{WBI}
	\frac{1}{2}\mathcal{L}_f |\nabla u|^2 &\ge& \frac{1}{m}(\mathcal{L}_f u)^2 + \langle \nabla u,\nabla\mathcal{L}_f u \rangle + Ric_f^{m-n}(\nabla u,\nabla u).
	\end{eqnarray}
\end{lemma}
We remark that the inequality (\ref{WBI}) is equivalent to the curvature dimension condition $CD(K,m)$ with respect to the operator $\mathcal{L}_f$ under the lower bound assumption $Ric_f^{m-n} \ge Kg$,  where $K\in \mathbb{R}$, $n <m \le \infty$ embedded in the next given inequality
\begin{align}\label{e2a}
\frac{1}{2}\mathcal{L}_f(|\nabla u|^2)- \langle \nabla u, \nabla \mathcal{L}_f u \rangle \ge \frac{1}{m}(\mathcal{L}_f u)^2  +  K|\nabla u|^2.
\end{align}
The assumption $Ric_f\ge Kg$ is clearly seen to imply $CD(K,\infty)$ by which it is obvious that \eqref{e2a} implies
\begin{align}\label{e3a}
\frac{1}{2}\mathcal{L}_f(|\nabla u|^2)- \langle \nabla u, \nabla \mathcal{L}_f u \rangle \ge  K|\nabla u|^2
\end{align}
but on the contrary  \eqref{e3a} does not yield  \eqref{e2a}. Thus, the lower bound condition on $Ric_f^{m-n}$ implies corresponding condition on $Ric_f$ but not the other way round.  Both the inequalities (\ref{e2a}) and (\ref{e3a}) are very useful in the analysis of diffusion operators and their geometry (see the books \cite{BGL,FWang}  and some references cited therein). The equation $Ric_f=\lambda g$ for some constant $\lambda$ defines the gradient Ricci solitons equation which has basic applications in the singularity analysis of the Ricci flow \cite{Hamilton-1,Perelman-1}, whilst the relation $Ric_f^{m-n}=K g$ corresponds to the quasi-Einstein equation which consists gradient Ricci solitons for $m=\infty$, and Einstein metrics  when $f$ is constant.  Quasi-Einstein metrics have been studied by several authors see \cite{CASE-1,CASE-2} for examples.

The next lemma is very important in the proof of Hessian estimate.
\begin{lemma}[\cite{Wang-hessian}]\label{lemma_hessian_prereq}
	For any two smooth function $u,v$ on $M$ we have 
	\begin{enumerate}
		\item $(\mathcal{L}_f-\partial_t)\frac{u}{v} = \frac{1}{v}(\mathcal{L}_f-\partial_t)f-\frac{u}{v^2}(\mathcal{L}_f-\partial_t)v-\frac{2}{v}\langle \nabla\frac{u}{v},\nabla v \rangle$,
		\item  $\mathcal{L}_f |\nabla^2u|^2 = 2|\nabla^2 u|\mathcal{L}_f|\nabla^2 u|+2|\nabla|\nabla^2 u||^2$.
	\end{enumerate}
\end{lemma}
\section{Main Results}
This section is divided into three subsections and each subsection contains a specific result related to the estimation. The Li-Yau type gradient estimates and its applications are discussed in subsection \ref{sec1}, Hamilton type gradient estimates and its application are discussed in subsection \ref{sec2}. Hessian estimate is presented in subsection \ref{sec3}.

\subsection{Li-Yau type estimation}\label{sec1}
For any $R>0$ and $p\in M$ we denote the geodesic ball at $p$ of radius $R$ by $\displaystyle B_p(R):=\{x\in M:d(x,p)\le R\}$ where $d(x,p)$ is the geodesic distance between $x,p$. The completeness of the manifold will assure the existence of geodesics.
\begin{lemma}\label{lemma1}
	Let $(M,g)$ be a complete $n$-dimensional Riemannian manifold with $Ric^{m-n}_f\ge -Kg$, where $K$ is a non-negative real number. If $u$ is a solution of (\ref{eq_heat_new}) then for any $\alpha \ge 1$ the quantity $\displaystyle F:=t(|\nabla u|^2-\alpha u_t -\alpha q-\alpha\tilde{G})$ satisfies
	\begin{eqnarray}
		\nonumber\label{eq_lemma1}\nonumber (\mathcal{L}_f - \partial_t)F &\ge& \frac{2t}{m}(\mathcal{L}_f u)^2 -2Kt|\nabla u|^2 - 2\langle \nabla F,\nabla u \rangle -\frac{F}{t}-\alpha t(\mathcal{L}_f q\\
		&&+\mathcal{L}_f \tilde{G})-2t(\alpha-1)\langle \nabla q,\nabla u\rangle-2t(\alpha-1)\langle \nabla \tilde{G},\nabla u\rangle.
	\end{eqnarray}
\end{lemma}
\begin{proof}
	Using weighted Bochner formula we get
	\begin{eqnarray}
		\label{eq_Lf}\nonumber \mathcal{L}_f F &\ge& \frac{2t}{m}(\mathcal{L}_f u)^2 + 2tRic_f^{m-n}(\nabla u,\nabla u)+ 2t\langle \nabla\mathcal{L}_f u,\nabla u \rangle - \alpha t\mathcal{L}_f u_t \\
		\nonumber &&- \alpha t\mathcal{L}_f q - \alpha t\mathcal{L}_f\tilde{G}\\
		\nonumber&\ge& \frac{2t}{m}(\mathcal{L}_fu)^2 - 2Kt |\nabla u|^2 -2\langle \nabla F,\nabla u\rangle -2t(\alpha-1)\langle \nabla u_t,\nabla u \rangle\\
		\nonumber && -2t(\alpha-1)\langle \nabla q,\nabla u\rangle -2t(\alpha-1)\langle \nabla\tilde{G},\nabla u \rangle -\alpha t \mathcal{L}_f u_t \\
		&&-\alpha t\mathcal{L}_f q - \alpha t\mathcal{L}_f \tilde{G}.
	\end{eqnarray}
	The last inequality is obtained by applying the Bakry-\'Emery Ricci curvature bound. Differentiating $F$ with respect to $t$ we get
	\begin{eqnarray}
	\label{eq_Ft} \nonumber F_t &=& |\nabla u|^2 - \alpha(u_t+q+\tilde{G})+t(\partial_t |\nabla u|^2 - \alpha(u_{tt}+q_t+\tilde{G}_t))\\
	&=& \frac{F}{t}+2t\langle\nabla u_t,\nabla u\rangle-\alpha t(u_{tt}+q_t+\tilde{G}_t)
	\end{eqnarray}
	Combining (\ref{eq_Lf}) and (\ref{eq_Ft}) we get 
	\begin{eqnarray}
		\nonumber(\mathcal{L}_f-\partial_t)F &\ge& \frac{2t}{m}(\mathcal{L}_fu)^2 - 2Kt |\nabla u|^2 -2\langle \nabla F,\nabla u\rangle -2t(\alpha-1)\langle \nabla u_t,\nabla u \rangle\\
		\nonumber && -2t(\alpha-1)\langle \nabla q,\nabla u\rangle -2t(\alpha-1)\langle \nabla\tilde{G},\nabla u \rangle \\
		\nonumber && -\alpha t \mathcal{L}_f u_t -\alpha t\mathcal{L}_f q - \alpha t\mathcal{L}_f \tilde{G}-F_t\\
		\nonumber&=&\frac{2t}{m}(\mathcal{L}_fu)^2 - 2Kt |\nabla u|^2 -2\langle \nabla F,\nabla u\rangle -2t(\alpha-1)\langle \nabla u_t,\nabla u \rangle\\
		\nonumber && -2t(\alpha-1)\langle \nabla q,\nabla u\rangle -2t(\alpha-1)\langle \nabla\tilde{G},\nabla u \rangle -\alpha t \mathcal{L}_f u_t -\alpha t\mathcal{L}_f q \\
		\nonumber && - \alpha t\mathcal{L}_f \tilde{G} -\frac{F}{t}-2t\langle\nabla u_t,\nabla u\rangle+\alpha t(u_{tt}+q_t+\tilde{G}_t)\\
		\nonumber&=&\frac{2t}{m}(\mathcal{L}_fu)^2 - 2Kt |\nabla u|^2 -2\langle \nabla F,\nabla u\rangle -2t(\alpha-1)\langle \nabla q,\nabla u\rangle \\
		\nonumber && -2t(\alpha-1)\langle \nabla\tilde{G},\nabla u \rangle -\alpha t(\mathcal{L}_f q + \mathcal{L}_f \tilde{G})-\frac{F}{t}-\alpha t \mathcal{L}_f u_t\\
   	 \label{neq1} && -2\alpha t\langle\nabla u_t,\nabla u\rangle+\alpha t(u_{tt}+q_t+\tilde{G}_t).
	\end{eqnarray}Note that 
\begin{eqnarray}
	\label{neq2}\nonumber u_{tt}+q_t+\tilde{G}_t &=& \frac{\partial}{\partial t}(\mathcal{L}_f u+|\nabla u|^2 - q -\tilde{G}) + q_t + \tilde{G}_t\\
	&=& \mathcal{L}_f u_t + 2\langle \nabla u_t,\nabla u \rangle.
\end{eqnarray} Hence we have
\begin{eqnarray}
	\label{neq3} -\alpha t \mathcal{L}_f u_t + \alpha t( u_{tt}+q_t+\tilde{G}_t) &=& 2\alpha t\langle \nabla u_t,\nabla u \rangle.
\end{eqnarray}Finally substituting (\ref{neq3}) into (\ref{neq1}) makes its last three terms vanish and thus we arrived at the desired result (\ref{eq_lemma1}).
\end{proof}
\begin{theorem}
	Let $(M,g)$ be a compact $n$-dimensional Riemannian manifold with $Ric_f^{m-n}\ge -K g$. Suppose that the boundary $\partial M$ of $M$ is convex whenever $\partial M\ne\emptyset$. If $w=e^u$ is a solution of $(\mathcal{L}_f-\partial_t)w =G(w)$ on $M\times (0,T]$, with Neumann boundary condition $\frac{\partial w}{\partial \nu} = 0$ on $\partial M \times (0,T]$, then
		\begin{eqnarray}
			\label{eqNTh1} \frac{|\nabla w|^2}{w^2}-\frac{2w_t}{w}-\frac{2G(w)}{w} &\le& \frac{m}{t} + \sqrt{2m\left|\mathcal{L}_f\left(\frac{G(w)}{w}\right)\right|}+mK,
		\end{eqnarray}
\end{theorem}
\begin{proof}
	For $\alpha=1$ and $q=0$ we have from Lemma \ref{lemma1}
	\begin{eqnarray}
	\nonumber (\mathcal{L}_f - \partial_t)F &\ge& \frac{2t}{m}(\mathcal{L}_f u)^2 -2Kt|\nabla u|^2 - 2\langle \nabla F,\nabla u \rangle - t\mathcal{L}_f \tilde{G} -\frac{F}{t},
	\end{eqnarray}
	where $F=t(|\nabla u|^2 - u_t - \tilde{G})=-t\mathcal{L}_f u$. Thus using $|x|\ge x$, for $x\in \mathbb{R}$ we get
	\begin{eqnarray}\label{eqN1}
	\nonumber (\mathcal{L}_f - \partial_t)F &\ge& \frac{2F^2}{mt} -2Kt|\nabla u|^2 - 2\langle \nabla F,\nabla u \rangle - t|\mathcal{L}_f \tilde{G}| -\frac{F}{t}\\
	 &=& -2\langle \nabla F,\nabla u \rangle+\frac{2F}{mt}\left(F-\frac{m}{2}-\frac{mt^2}{2F}|\mathcal{L}_f \tilde{G}|\right)-2Kt|\nabla u|^2.
	\end{eqnarray}
	Let us suppose that $\displaystyle \mathcal{L}_f\tilde{G}\ge 0$. We claim that 
	\begin{eqnarray}\label{eqClaim}
		\frac{2F}{mt}\left(F-\frac{m}{2}-\frac{mt^2}{2F}|\mathcal{L}_f \tilde{G}|\right)-2Kt|\nabla u|^2\le 0.
	\end{eqnarray} Assume the contrary, let $(x_0,t_0)\in M\times (0,T]$ be a point in which 
	\begin{eqnarray}\label{eqN2}
	 \frac{2F(x_0,t_0)}{mt_0}\left(F(x_0,t_0)-\frac{m}{2}-\frac{mt_0^2}{2F(x_0,t_0)}|\mathcal{L}_f \tilde{G}|\right)-2Kt_0|\nabla u|^2\le 0.
	\end{eqnarray}
	Now there are two cases. If $x_0$ is an interior point of $M$, then 
	\begin{eqnarray}
	\nonumber (\mathcal{L}_f-\partial_t) F\Big|_{(x_0,t_0)}\le 0 = \nabla F\Big|_{(x_0,t_0)}.
	\end{eqnarray} 
	Hence from (\ref{eqN1}) at $(x_0,t_0)$ we deduce
	\begin{eqnarray}
	\nonumber 0\ge (\mathcal{L}_f-\partial_t)F&\ge&\frac{2F(x_0,t_0)}{mt_0}\left(F(x_0,t_0)-\frac{m}{2}-\frac{mt_0^2}{2F(x_0,t_0)}|\mathcal{L}_f \tilde{G}|\right)-2Kt_0|\nabla u|^2\\
	\nonumber &>& 0,
	\end{eqnarray}
	which is a contradiction to (\ref{eqN2}). Thus (\ref{eqClaim}) holds for any interior point $x_0\in M$.\\
	Next we suppose that $x_0\in\partial M$. Following \cite{Yi-Li}, by strong maximum principle we have $\displaystyle \frac{\partial F}{\partial \nu}(x_0,t_0)>0$. Let $\{e_i:i=1,2,\cdots,n\}$ be an orthonormal basis of $TM$, where $e_m :=\frac{\partial}{\partial \nu}$. Thus 
	\begin{eqnarray}
	\frac{\partial F}{\partial \nu} &=& 2t\sum_{j=1}^{n-1}u_j \frac{\partial u_j}{\partial \nu} + 2t\frac{\partial u}{\partial \nu}\frac{\partial^2 u}{\partial \nu^2} - \frac{\partial u_t}{\partial \nu} - \frac{\partial \tilde{G}}{\partial \nu}.
	\end{eqnarray}
	Under Neumann boundary condition $\frac{\partial w}{\partial \nu} = 0$ on $\partial M$ which implies $\frac{\partial u}{\partial \nu}=0$ on $\partial M$ and consequently
	\begin{eqnarray}
	\nonumber \frac{\partial F}{\partial \nu} &=& 2t\sum_{j=1}^{n-1}u_j \frac{\partial u_j}{\partial \nu}\\
	\nonumber &=& -2t\sum_{j=1}^{n-1}\sum_{k=1}^{n-1} h_{jk}u_ju_k\\
	\nonumber &=& -2t\ {\rm II}(\nabla u,\nabla u).
	\end{eqnarray}
	The second line was due to $\displaystyle u_{j\nu} = -\sum_{1\le k\le n-1}h_{jk}u_k$, where $h_{jk}$ are components of the second fundamental form ${\rm II}$ of $\partial M$. It follows immediately that ${\rm II}(\nabla u,\nabla u)(x_0,t_0)<0$, contradicting the convexity of $\partial M$. Hence (\ref{eqClaim}) holds on whole of $M$ under the obvious restriction $|\mathcal{L}_f\tilde{G}|\ge 0$. So in this case we have a quadratic in $F$ as
	\begin{eqnarray}
	F^2 - \frac{m}{2}F-mt^2\left(\frac{1}{2}|\mathcal{L}_f\tilde{G}|+K|\nabla u|^2\right) &\le& 0.
	\end{eqnarray}
	For $a,b,x\ge 0$ and $x^2\le ax+b$ we have $x\le\frac{a}{2}+\sqrt{b+\left(\frac{a}{2}\right)^2}\le a+\sqrt{b}$. Treating $F$ as $x$ we get
	\begin{eqnarray}
	\nonumber F &\le& \frac{m}{2}+t\sqrt{\frac{m}{2}|\mathcal{L}_f \tilde{G}|+mK|\nabla u|^2}\\
	&\le& \frac{m}{2}+t\sqrt{\frac{m}{2}|\mathcal{L}_f \tilde{G}|}+t\sqrt{mK}|\nabla u|.
	\end{eqnarray}
	By Young's inequality and from the definition of $F$ we find 
	\begin{eqnarray}
		\nonumber |\nabla u|^2-u_t-\tilde{G} &\le& \frac{m}{2t}+\sqrt{\frac{m}{2}|\mathcal{L}_f \tilde{G}|}+\frac{mK}{2}+\frac{|\nabla u|^2}{2},\\
		\text{or, } \frac{1}{2}|\nabla u|^2-u_t-\tilde{G} &\le& \frac{m}{2t}+\sqrt{\frac{m}{2}|\mathcal{L}_f \tilde{G}|}+\frac{mK}{2},
	\end{eqnarray}
	which is the same as (\ref{eqNTh1}) after substituting $u=\log w$.
\end{proof}
\begin{remark}
	Taking $K=0$ will produce the classical estimates as showed in \cite{Li-Yau,Yi-Li}.
\end{remark}

\begin{theorem}\label{theorem1}
	If $u$ is a positive solution of the equation (\ref{eq_heat_new}) on a complete weighted Riemannian manifold $M$ with boundary and $Ric_f^{m-n}\ge -Kg$, on $B_p(2R)$, for some non-negative constant $K$, then
	\begin{eqnarray}\label{theorem1_eq1}
	\nonumber && |\nabla u|^2 - \alpha u_t -\alpha q - \alpha\tilde{G}(u)\\
	 \nonumber &\le& \Bigg[\frac{m\alpha^2A}{2(1-\epsilon)}+\frac{m\alpha^2K}{(1-\epsilon)(\alpha-1)}+\frac{m^2\alpha^4C_1^2}{4\epsilon R^2(1-\epsilon)(\alpha-1)}+\frac{m\alpha^2}{2(1-\epsilon)}\\
	\nonumber &&+\left(\frac{m\alpha^3(\mathcal{L}_f q+\mathcal{L}_f\tilde{G})}{2(1-\epsilon)}+\frac{\alpha^2(\alpha-1)m}{1-\epsilon}(|\nabla q|^2+|\nabla \tilde{G}|^2)\right)^\frac12\Bigg]\\
	&&+\frac{m\alpha^2}{2t(1-\epsilon)}
	\end{eqnarray}
on $B_p(2R)$ for any $0<\epsilon<1$, where $\displaystyle A=\frac{(m-1)C_1 (1+R\sqrt{K})+C_2+2C_1^2}{R^2}$.
\end{theorem}
\begin{proof}
	Following \cite{Chen-Chen,Li-Yau,Ma,Zhu-Li}, we consider a cut-off function $\psi:[0,\infty]\to \mathbb{R}$ defined by
	\begin{equation}
	\nonumber \psi(r) = \begin{cases}
	1,\ r\in[0,1],\\
	0,\ r\in[2,\infty),
	\end{cases}
	\end{equation}
	with $\displaystyle -C_1\le\frac{\psi'(r)}{\sqrt{\psi(t)}}\le 0$ and $\displaystyle \psi(r)\ge -C_2$, for some positive constant $C_1,C_2$. Define $\displaystyle \varphi(x):=\psi\left(\frac{d(x,p)}{R}\right)$. Calabi's trick \cite{CALABI,Cheng-Yau} enables us to assume the smoothness of $\varphi$ in $B_p(2R)$. By \cite[Corollary~3.3]{Yi-Li} we find $\frac{|\nabla \varphi|^2}{\varphi}\le \frac{C_1^2}{R^2}$ and $\mathcal{L}_f \varphi \ge -\frac{(m-1)C_1(1+R\sqrt{K})+C_2}{R^2}$. Fix a time $T_1\le T$ and let $(x_0,t_0)\in M\times [0,T_1]$ where $\mathcal{G}=\varphi F$ achieves its maximum. Hence at $(x_0,t_0)$ we have the following
	\begin{equation}
	\nonumber \nabla \mathcal{G}=0,\ \ (\mathcal{L}_f-\partial_t)\mathcal{G}\le 0.
	\end{equation}
	We see that $\nabla\mathcal{G}=0\implies \nabla F=-\frac{F}{\varphi}\nabla\varphi$. Using Lemma \ref{lemma1} we infer
	\begin{eqnarray}\label{eq1}
	\nonumber\mathcal{L}_f\mathcal{G} &=& F\mathcal{L}_f\varphi +2\langle\nabla\varphi,\nabla F\rangle+\varphi\mathcal{L}_fF\\
	\nonumber &\ge& -AF+\varphi\mathcal{L}_f F\\
	\nonumber (\mathcal{L}_f-\partial_t)\mathcal{G} &\ge& -AF +\frac{2t_0}{m} \varphi \nonumber(\mathcal{L}_f u)^2 - 2Kt_0\varphi |\nabla u|^2 -2F \langle \nabla\varphi,\nabla u \rangle \\
	\nonumber&&-\alpha t_0\varphi (\mathcal{L}_fq+\mathcal{L}_f\tilde{G})-\frac{\varphi F}{t_0} -2t_0 \varphi (\alpha-1)\langle \nabla q,\nabla u \rangle\\
	&& -2t_0 \varphi (\alpha-1)\langle \nabla \tilde{G},\nabla u \rangle.
	\end{eqnarray}
	Since $(\mathcal{L}_f-\partial_t)\mathcal{G}\le 0$ thus the above equation reduces to 
	\begin{eqnarray}\label{eq2}
	\nonumber 0 &\ge& -AF +\frac{2t_0}{m} \varphi \nonumber(\mathcal{L}_f u)^2 - 2Kt_0\varphi |\nabla u|^2 -2F \langle \nabla\varphi,\nabla u \rangle -\alpha t_0\varphi (\mathcal{L}_fq+\mathcal{L}_f\tilde{G})\\
	&&-\frac{\varphi F}{t_0}-2t_0 \varphi (\alpha-1)\langle \nabla q,\nabla u \rangle -2t_0 \varphi (\alpha-1)\langle \nabla \tilde{G},\nabla u \rangle. 
	\end{eqnarray}
	Following \cite{Chen-Chen,YANG,Zhu-Li} we set $$\mu=\frac{|\nabla u|^2}{F}\Big|_{(x_0,t_0)}\ge 0.$$ Multiplying (\ref{eq2}) with $\varphi t_0$ and using the above relation we get
	\begin{eqnarray}\label{eq3}
	\nonumber A\mathcal{G}t_0 &\ge& \frac{2t_0^2}{m}\varphi^2(\mathcal{L}_f u)^2 - 2Kt_0^2\varphi \mu \mathcal{G} + 2\mathcal{G}t_0 \langle\nabla\varphi,\nabla u\rangle - \alpha t_0^2\varphi^2(\mathcal{L}_fq+\mathcal{L}_f\tilde{G})-\varphi \mathcal{G}\\
	\nonumber&& -2t_0^2 \varphi^2 (\alpha-1)\langle \nabla q,\nabla u \rangle -2t_0^2 \varphi^2 (\alpha-1)\langle \nabla \tilde{G},\nabla u \rangle.
	\end{eqnarray}
	By Cauchy-Schwarz inequality
	\begin{eqnarray}
	\label{eqCS1} \langle\nabla u,\nabla\varphi\rangle &\le& |\nabla u||\nabla\varphi|\ \le\  \frac{C_1}{R}\mu^\frac12 \mathcal{G}^\frac12,\\
	\label{eqCS2} \langle \nabla q,\nabla u \rangle &\le& |\nabla q| |\nabla u|\ \le\ |\nabla q|^2 + \frac{|\nabla u|^2}{4},\\
	\label{eqCS3} \langle \nabla \tilde{G},\nabla u \rangle &\le& |\nabla \tilde{G}| |\nabla u|\ \le\ |\nabla \tilde{G}|^2 + \frac{|\nabla u|^2}{4}.
	\end{eqnarray}
	Using (\ref{eqCS1}), (\ref{eqCS2}) and (\ref{eqCS3}) in (\ref{eq3}) we find that
	\begin{eqnarray}\label{eq4_old}
	\nonumber A\mathcal{G}t_0 &\ge& \frac{2t_0^2}{m}\varphi^2(\mathcal{L}_f u)^2 - (2K+\alpha-1)t_0^2\varphi \mu \mathcal{G} - \frac{2C_1}{R}\mu^\frac12 t_0 \mathcal{G}^\frac32 - \alpha t_0^2\varphi^2(\mathcal{L}_fq+\mathcal{L}_f\tilde{G})\\
	&&-\varphi \mathcal{G} - 2t_0^2 \varphi^2 (\alpha-1) (|\nabla q|^2+|\nabla \tilde{G}|^2).
	\end{eqnarray}
	In terms of $F$ equation (\ref{eq_heat_new}) becomes
	\begin{eqnarray}\label{eq5}
	\nonumber (\mathcal{L}_f u)^2 &=& (u_t-|\nabla u|^2+q+\tilde{G})^2\\
	&=& F^2\left(\mu - \frac{\mu t_0-1}{\alpha t_0}\right)^2
	\end{eqnarray}
	Using (\ref{eq5}) in (\ref{eq4_old})  we deduce
	\begin{eqnarray}\label{eq6}
	\nonumber A\mathcal{G}t_0 &\ge& \frac{2t_0^2}{m}(\mu-\frac{\mu t_0-1}{\alpha t_0})^2\mathcal{G}^2 - (2K+\alpha-1)t_0^2\varphi \mu \mathcal{G} - \frac{2C_1}{R}\mu^\frac12 t_0 \mathcal{G}^\frac32 \\
	&&- \alpha t_0^2\varphi^2(\mathcal{L}_f q+ \mathcal{L}_f \tilde{G})-\varphi \mathcal{G}- 2t_0^2 \varphi^2 (\alpha-1) (|\nabla q|^2+ |\nabla \tilde{G}|^2).
	\end{eqnarray}
	For $0<\epsilon<1$ we have by Young's inequality
	\begin{eqnarray}\label{eq7}
	\frac{2C_1 t_0}{R}\mu^\frac12 \mathcal{G}^\frac32 &\le& \frac{2\epsilon}{m\alpha^2}[1+(\alpha-1)\mu t_0]^2\mathcal{G}^2 + \frac{m\alpha^2C_1^2 t_0^2 \mu}{2\epsilon R^2 [1+(\alpha-1)\mu t_0]^2}\mathcal{G}.
	\end{eqnarray} 
	Using (\ref{eq7}) in (\ref{eq6}) we deduce
	\begin{eqnarray}\label{eq8}
	\nonumber && \left[\frac{2(1-\epsilon)}{m\alpha^2}(1+(\alpha-1)\mu t_0)^2\right]\mathcal{G}^2\\
	\nonumber &\le& \left[At_0+\varphi+(2K+\alpha-1)t_0^2\varphi\mu+\frac{m\alpha^2C_1^2 t_0^2\mu}{2\epsilon R^2[1+(\alpha-1)\mu t_0]^2}\right]\mathcal{G}\\
	&&+\alpha t_0^2\varphi(\mathcal{L}_f q + \mathcal{L}_f \tilde{G})+ 2t_0^2 \varphi^2 (\alpha-1) (|\nabla q|^2+|\nabla \tilde{G}|^2).
	\end{eqnarray}
	By definition of $\varphi$ we see $0\le\varphi\le 1$ and $1+(\alpha-1)\mu t_0\ge 1$. Hence (\ref{eq8}) further reduces to
	\begin{eqnarray}\label{eq9}
	\nonumber \frac{2(1-\epsilon)}{m\alpha^2}\mathcal{G}^2 &\le& \left[At_0+1+t_0+\frac{2Kt_0}{\alpha-1}+\frac{m\alpha^2C_1^2 t_0}{2\epsilon R^2(\alpha-1)}\right]\mathcal{G}+\alpha t_0^2( \mathcal{L}_f q + \mathcal{L}_f \tilde{G})\\
	&&+ 2t_0^2 (\alpha-1) ( |\nabla q|^2 + |\nabla \tilde{G}|^2),
	\end{eqnarray}
	at $(x_0,t_0)$. Multiplying both sides with $\frac{m\alpha^2}{2(1-\epsilon)}$ and rearranging the terms we get
	
	\begin{eqnarray}
	\nonumber \mathcal{G}^2 &\le& \left[\frac{m\alpha^2At_0}{2(1-\epsilon)}+\frac{m\alpha^2(1+t_0)}{2(1-\epsilon)}+\frac{m\alpha^2Kt_0}{(1-\epsilon)(\alpha-1)}+\frac{m^2\alpha^4C_1^2t_0}{4\epsilon R^2(1-\epsilon)(\alpha-1)}\right]\mathcal{G}\\
	&&+\frac{m\alpha^3 t_0^2}{2\epsilon R^2(1-\epsilon)}(\mathcal{L}_f q+\mathcal{L}_f \tilde{G})+\frac{\alpha^2(\alpha-1)mt_0^2}{1-\epsilon}(|\nabla q|^2 + |\nabla \tilde{G}|^2).
	\end{eqnarray}
	For $a,b,x\ge 0$ whenever $x^2\le ax+b$ we have $x\le \frac{a}{2}+\sqrt{\frac{b}{2}+\left(\frac{a}{2}\right)^2}\le \frac{a}{2}+\sqrt{b}+\frac{a}{2}\le a+\sqrt{b}$. Hence from the above quadratic in $\mathcal{G}$ we get
	\begin{eqnarray}
	\nonumber \mathcal{G} &\le& \frac{m\alpha^2}{2(1-\epsilon)}+ \left[\frac{m\alpha^2A}{2(1-\epsilon)}+\frac{m\alpha^2K}{(1-\epsilon)(\alpha-1)}+\frac{m^2\alpha^4C_1^2}{4\epsilon R^2(1-\epsilon)(\alpha-1)}+\frac{m\alpha^2}{2(1-\epsilon)}\right]t_0\\
	&&+\left(\frac{m\alpha^3(\mathcal{L}_f q+\mathcal{L}_f\tilde{G})}{2(1-\epsilon)}+\frac{\alpha^2(\alpha-1)m}{1-\epsilon}(|\nabla q|^2+|\nabla\tilde{G}|^2)\right)^\frac12 t_0
	\end{eqnarray}
	at $(x_0,t_0)$. By construction of $\varphi$, $F\le \mathcal{G}(x_0,t_0)$ on $B_p(R)\times [0,T']$. Since $T'$ is arbitrary so we have
	\begin{eqnarray}\label{eqNew1}
	\nonumber F &\le& \Bigg[\frac{m\alpha^2A}{2(1-\epsilon)}+\frac{m\alpha^2K}{(1-\epsilon)(\alpha-1)}+\frac{m^2\alpha^4C_1^2}{4\epsilon R^2(1-\epsilon)(\alpha-1)}+\frac{m\alpha^2}{2(1-\epsilon)}\\
	\nonumber &&+\left(\frac{m\alpha^3(\mathcal{L}_fq+\mathcal{L}_f\tilde{G})}{2(1-\epsilon)}+\frac{\alpha^2(\alpha-1)m}{1-\epsilon}(|\nabla q|^2+|\nabla\tilde{G}|^2)\right)^\frac12\Bigg] t\\
	&& +\frac{m\alpha^2}{2(1-\epsilon)}.
	\end{eqnarray}
	Putting the value of $F$ in \eqref{eqNew1} completes the proof.
\end{proof}
\noindent As an immediate consequence we have the results.
\begin{corollary}\label{corr1}
	If $u$ is a positive solution of the equation (\ref{eq_heat_new}) on a complete weighted Riemannian manifold $M$ with boundary and $Ric_f^{m-n}\ge -Kg$ on $M\times (0,T]$ then
	\begin{eqnarray}\label{corr1_eq1}
	\nonumber && |\nabla u|^2 - \alpha u_t -\alpha q - \alpha\tilde{G}(u)\\
	\nonumber &\le& \frac{m\alpha^2}{2t(1-\epsilon)} + \Bigg[\frac{m\alpha^2K}{(1-\epsilon)(\alpha-1)}+\frac{m\alpha^2}{2(1-\epsilon)}\\
	&&+\left(\frac{m\alpha^3(\mathcal{L}_f q+\mathcal{L}_f \tilde{G})}{2(1-\epsilon)}+\frac{\alpha^2(\alpha-1)m}{1-\epsilon}(|\nabla q|^2+|\nabla\tilde{G}|^2)\right)^\frac12\Bigg],
	\end{eqnarray}
	on $M\times (0,T]$ for any $0<\epsilon<1$.
\end{corollary}
\begin{proof}
	Letting $R\to+\infty$ in Theorem \ref{theorem1} we get the above result.
\end{proof}
For the next result we consider a space-time curve $\gamma:(0,T]\to M$ joining the points $(x_1,t_1),\ (x_2,t_2)\in M\times (0,T]$ where $0<t_1<t_2<T$ such that $\gamma(t_i)=x_i$, $i=1,2$. To derive the Harnack type inequality and compare heat between two different points on the manifold let us define for $\alpha>1$ and $0<\epsilon<1$,
\begin{eqnarray}\label{eq_lambda}
\nonumber \Lambda_{\alpha,\epsilon} &=& \left(\frac{m\alpha^3(\mathcal{L}_f q+\mathcal{L}_f \tilde{G})}{2(1-\epsilon)}+\frac{\alpha^2(\alpha-1)m}{1-\epsilon}(|\nabla q|^2+|\nabla\tilde{G}|^2)\right)^\frac12 +\frac{m\alpha^2K}{(1-\epsilon)(\alpha-1)}\\
&&+\frac{m\alpha^2}{2(1-\epsilon)},
\end{eqnarray}
The idea behind such consideration is to make our results more compact by shortening its length.
\begin{corollary}
	If $w=e^u$ is a solution of (\ref{eq_heat_old}) on a complete weighted Riemannian manifold $M$ then for any two points $(x_1,t_1),(x_2,t_2)\in M\times (0,T]$ we have the Harnack type inequality for $f$ given by
	\begin{eqnarray}
	\nonumber \frac{w(x_1,t_1)}{w(x_2,t_2)} &\le&\left(\frac{t_1}{t_2}\right)^{\frac{m\alpha}{2(1-\epsilon)}}\exp\Big\{\inf_{\gamma}\int_{t_1}^{t_2}\Big(\frac{\Lambda_{\alpha,\epsilon}}{\alpha}+\frac{\alpha-2}{2\alpha}|\nabla u|^2\\
	&& +\frac{1}{2}|\dot{\gamma}(t)|^2-q(x,t)-G(w)\Big)dt\Big\},
	\end{eqnarray} 
	where the infimum is taken over all curves joining $(x_1,t_1),(x_2,t_2)$.
\end{corollary}
\begin{proof}
	Combining all the cases of Corollary \ref{corr1} gives 
	\begin{eqnarray}
	\nonumber |\nabla u|^2-\alpha u_t -\alpha q - \alpha\tilde{G}(u) &\le& \Lambda_{\alpha,\epsilon} + \frac{m\alpha^2}{2t(1-\epsilon)}.
	\end{eqnarray}
	Upon rearrangement we have
	\begin{eqnarray}\label{eq10}
	\frac{\partial u}{\partial t}&\ge& \frac{|\nabla u|^2}{\alpha}+q+\tilde{G}-\frac{\Lambda_{\alpha,\epsilon}}{\alpha}-\frac{m\alpha}{2t(1-\epsilon)}.
	\end{eqnarray}
	On the other hand
	\begin{eqnarray}\label{eq11}
	\nonumber \frac{du}{dt} &=& \frac{\partial u}{\partial t}+\langle\nabla u,\dot{\gamma}(t)\rangle\\
	\nonumber &\ge& \frac{\partial u}{\partial t}-\langle\nabla u,\dot{\gamma}(t)\rangle\\
	\nonumber &\ge& \frac{2-\alpha}{2\alpha}|\nabla u|^2+q+\tilde{G}-\frac{\Lambda_{\alpha,\epsilon}}{\alpha}-\frac{m\alpha}{2t(1-\epsilon)}-\frac{1}{2}|\dot{\gamma}(t)|^2.
	\end{eqnarray}
	Here we have used (\ref{eq10}) and an elementary inequality $\frac{1}{2}|\nabla u|^2-\langle\nabla u,\dot{\gamma}(t)\rangle\ge-\frac12 |\dot{\gamma}(t)|^2$. Integrating (\ref{eq11}) between $(x_1,t_1)$ and $(x_2,t_2)$ along the path $\gamma$ we get
	\begin{eqnarray}
	\nonumber u(x_2,t_2)-u(x_1,t_1) &\ge& \sup_{\gamma}\int_{t_1}^{t_2}\left(\frac{2-\alpha}{2\alpha}|\nabla u|^2+q+\tilde{G}-\frac{\Lambda_{\alpha,\epsilon}}{\alpha}-\frac{1}{2}|\dot{\gamma}(t)|^2\right)dt\\
	&&-\frac{m\alpha}{2(1-\epsilon)}\ln(t)\Big|_{t_1}^{t_2}.
	\end{eqnarray}
	Exponentiating both sides we find 
	\begin{eqnarray}
	\nonumber \frac{e^{u(x_2,t_2)}}{e^{u(x_1,t_1)}} &\ge& \exp\left\{\sup_{\gamma}\int_{t_1}^{t_2}\left(\frac{2-\alpha}{2\alpha}|\nabla u|^2+q+\tilde{G}-\frac{\Lambda_{\alpha,\epsilon}}{\alpha}-\frac{1}{2}|\dot{\gamma}(t)|^2\right)dt\right\}\left(\frac{t_2}{t_1}\right)^{-\frac{m\alpha}{2(1-\epsilon)}}.
	\end{eqnarray}
	Rewriting in terms of $w(=e^u)$ and reciprocating both sides completes the proof.
\end{proof}
\begin{remark}\label{remark_important}
	We made the above theorem flexible in the sense that no particular choice of $q$ and $\tilde{G}(u)$ are considered. In the introduction section we saw different values of $\tilde{G}$ and their importance in various fields of science. One can find Li-Yau type estimation and Harnack type inequalities for such cases by using the value of $\nabla\tilde{G},\ \mathcal{L}_f\tilde{G}$ and repeating the aforementioned process. But it should be mentioned that just by plugging the values of $\nabla\tilde{G},\ \mathcal{L}_f\tilde{G}$ in (\ref{theorem1_eq1}) will not produce any desired result. Because of the presence of terms like $|\nabla w|^2$, $\mathcal{L}_f w$ and constants. They need to be eliminated during the procedure with the help of Young's inequality and necessary restrictions (like assuming bounds for constants $a,p,q$ etc.). No matter what happens one must restraint himself from taking any assumption that involves bounds of $|\nabla u|$, because our objective is to find such a bound hence we can not assume that.
\end{remark}

\subsection{Hamilton type estimate and Liouville type theorem}\label{sec2}
Next we establish a Liouville type theorem for positive solution of $w$ of (\ref{eq_heat_old}). For this we need the following lemma.
\begin{lemma}\label{eqLTTlemma}
	For any smooth positive solution $w$ of (\ref{eq_heat_old}) on a compact weighted Riemannian manifold $M$ with $Ric_f\ge -Kg$, for some $K\ge 0$, the quantity $\frac{|\nabla w|^2}{w}$ satisfies
	\begin{eqnarray}\label{eqLTTlemmaeq}
		\nonumber (\partial_t-\mathcal{L}_f)\frac{|\nabla w|^2}{w} &\le& \frac{|\nabla w|^2}{w}\left( 2K + \frac{G(w)}{w} -  G'(w) + |q| + 1 \right) +|\nabla q|^2 w.\\
	\end{eqnarray}
\end{lemma}
\begin{proof}
	Differentiating $\frac{|\nabla w|^2}{w}$ with respect to $t$ we obtain
	\begin{eqnarray}
		\label{lem_eq_1} \nonumber \partial_t\left(\frac{|\nabla w|^2}{w}\right) &=& \frac{\partial_t |\nabla w|^2}{w}-\frac{|\nabla w|^2}{w^2}\partial_t w\\
		\nonumber &=& \frac{2}{w}(\langle \nabla w,\nabla\mathcal{L}_f w \rangle - \langle \nabla w,\nabla(qw+G(w))\rangle) \\
		&&- \frac{|\nabla w|^2}{w^2}(\mathcal{L}_f w - qw -G(w)).
	\end{eqnarray}
A consequence of direct calculation gives
\begin{eqnarray}
	\nonumber \mathcal{L}_f \frac{|\nabla w|^2}{w} &=& \frac{\mathcal{L}_f |\nabla w|^2}{w} + |\nabla w|^2 \mathcal{L}_f(\frac1w) + 2\langle \nabla|\nabla w|^2,\nabla(\frac1w) \rangle\\
	\nonumber &=& \frac{\mathcal{L}_f|\nabla w|^2}{w} + |\nabla w|^2 (-\frac{1}{w^2}\mathcal{L}_f w+\frac{2|\nabla w|^2}{w^3}) -\frac{4}{w^2}\text{Hess }w(\nabla w,\nabla w).
\end{eqnarray}
Applying weighted Bochner formula (Lemma \ref{lemma_Bochner}) we have
\begin{eqnarray}\label{lem_eq_2}
\nonumber \mathcal{L}_f \frac{|\nabla w|^2}{w} &=& \frac{2}{w} Ric_f(\nabla w,\nabla w) + \frac{2}{w}\left( |\text{Hess }w|^2 - \frac{2}{w}\text{Hess }w(\nabla w,\nabla w) + \frac{|\nabla w|^4}{w^2} \right) \\
&& - \frac{|\nabla w|^2}{w^2}\mathcal{L}_f w + \frac{2}{w}\langle \nabla\mathcal{L}_f w,\nabla w \rangle.
\end{eqnarray}
Combining (\ref{lem_eq_1}) and (\ref{lem_eq_2}) we arrive at
\begin{eqnarray}\label{eq_nabwbyw}
	\nonumber && (\partial_t-\mathcal{L}_f)\frac{|\nabla w|^2}{w}\\
	\nonumber &=& \frac{2}{w} \langle \nabla\mathcal{L}_f w,\nabla w \rangle - \frac{2}{w} \langle \nabla(qw+G(w)), \nabla w \rangle - \frac{|\nabla w|^2}{w^2}\mathcal{L}_f w\\
	\nonumber && + \frac{|\nabla w|^2}{w^2}(qw+G(w)) - \frac{2}{w} \left( |\text{Hess }w|^2 - \frac{2}{w}\text{Hess }w(\nabla w,\nabla w) + \frac{|\nabla w|^4}{w^2}\right)\\
	 &&- \frac{2}{w} \langle \nabla\mathcal{L}_f w,\nabla w \rangle -\frac{2}{w} Ric_f(\nabla w,\nabla w) + \frac{|\nabla w|^2}{w^2}\mathcal{L}_f w.
\end{eqnarray}
Note that $\nabla G(w) = G'(w)\nabla w$, where $G'(w) = \dfrac{\partial}{\partial w}G(w)$. Then one has $$\frac{2}{w}\langle \nabla G(w), \nabla w \rangle = \frac{2|\nabla w|^2}{w}G'(w).$$ Using the above, the curvature condition $Ric_f\ge -Kg$ and the inequality $$\frac{2}{w}|\text{Hess }w|^2-\frac{4}{w^2}\text{Hess } w(\nabla w,\nabla w)+\frac{2|\nabla w|^4}{w^3}=\frac{2}{w}\Big|\text{Hess }w - \frac{\nabla w\otimes \nabla w}{w}\Big|^2\ge 0,$$ into (\ref{eq_nabwbyw}) we have
\begin{equation}\label{eq_2.39}
(\partial_t-\mathcal{L}_f)\frac{|\nabla w|^2}{w} \le \frac{|\nabla w|^2}{w}\left( 2K + \frac{G(w)}{w} -  G'(w) + |q| \right) -2\langle \nabla q,\nabla w \rangle.
\end{equation}
Finally to simplify the last term we apply Cauchy--Schwarz and Young's inequality. Thus we get 
\begin{eqnarray}\label{eq_Youngcauchy}
\nonumber -2\langle \nabla q,\nabla w \rangle &\le& 2 |\nabla q||\nabla w|\\
\nonumber &=& 2 |\nabla q|\sqrt{w}\ \frac{|\nabla w|}{\sqrt{w}}\\
 &\le& |\nabla q|^2 w + \frac{|\nabla w|^2}{w}.
\end{eqnarray}
Using (\ref{eq_Youngcauchy}) in (\ref{eq_2.39}) we get the result (\ref{eqLTTlemmaeq}).
\end{proof}
\begin{theorem}\label{TheoremLvTh}
	Let $M$ be a compact weighted Riemannian manifold with $Ric_f\ge-Kg$ where $K\ge 0$. If $w$ is a positive solution of (\ref{eq_heat_old}) with $0<w\le \frac{A}{e}$, $|q|\le \theta_1,\ |G|\le\theta_2|w|$, $|G'(w)|\le\theta_3$, $|\nabla q|\le \theta_4$ for some $\theta_1,\theta_2,\theta_3,\theta_4\ge 0$, on $M\times (0,T]$ then 
	\begin{eqnarray}
	\label{eqTheoremLvTh}
	\frac{|\nabla w|^2}{w} &\le& \frac{A}{e} \left((\ln\frac{A}{w}-1)(\theta_1+\theta_2)+\ln(\frac{A}{w}) + \frac{\theta_4^2}{\xi}\right)\left(\frac1t + \xi\right),
	\end{eqnarray}
	where $\xi = 2K + \theta_1+\theta_2 + \theta_3 + 1$.
\end{theorem}
\begin{proof}
	First note that $0<w\le\frac{A}{e}\implies \ln\frac{A}{w}-1\ge 0$. By a direct calculation we find that
	\begin{eqnarray}
	\label{eq1_Thm1}\nonumber \partial_t (w\ln\frac{A}{w}) &=& \partial_tw \ \ln\frac{A}{w} + w\ \partial_t(\ln\frac{A}{w})\\
	&=& \partial_tw\ (\ln\frac{A}{w}-1).
	\end{eqnarray}
	Similarly
	\begin{eqnarray}\label{eq2_Thm1}
	\nonumber\mathcal{L}_f (w\ln\frac{A}{w}) &=& \mathcal{L}_f w\ \ln\frac{A}{w} + w\ \mathcal{L}_f\ln\frac{A}{w} + 2\langle \nabla w,\nabla\ln\frac{A}{w} \rangle\\
	&=& (\ln\frac{A}{w}-1)\mathcal{L}_fw-\frac{|\nabla w|^2}{w}.
	\end{eqnarray}
	Combining (\ref{eq1_Thm1}) and (\ref{eq2_Thm1}) gives
	\begin{eqnarray}
	\nonumber (\partial_t-\mathcal{L}_f)(w\ln\frac{A}{w}) &=& (\ln\frac{A}{w}-1)(\partial_t-\mathcal{L}_f)w +\frac{|\nabla w|^2}{w}\\
	&=& (\ln\frac{A}{w}-1)(-qw-G(w)) +\frac{|\nabla w|^2}{w}.
	\end{eqnarray}
	Consider the function
	\begin{eqnarray}
	\nonumber F:= \phi \frac{|\nabla w|^2}{w} - w\ln\frac{A}{w},
	\end{eqnarray}
	where $\phi$ is a time dependent function with $\phi(0) = 0$. Using Lemma \ref{eqLTTlemma} we see that
	\begin{eqnarray}\label{eq3_Theorem1}
	\nonumber	(\partial_t-\mathcal{L}_f)F &=& \phi_t\frac{|\nabla w|^2}{w} + \phi\  \partial_t\frac{|\nabla w|^2}{w} - \partial_t(w\ln\frac{A}{w}) - \phi \mathcal{L}_f(\frac{|\nabla w|^2}{w}) + \mathcal{L}_f(w \ln\frac{A}{w})\\
	\nonumber &\le& (\phi_t + \phi (2K+\frac{G(w)}{w}-G'(w)+|q|+1)-1)\frac{|\nabla w|^2}{w}\\
	&& +(\ln\frac{A}{w}-1)(qw+G(w)) +\phi |\nabla q|^2 w.
	\end{eqnarray}
	Using the bounds of $q,\ G(w)$ and $G'(w)$ we obtain
	\begin{eqnarray}
	\nonumber (\partial_t-\mathcal{L}_f)F &\le& (\phi_t + \phi (2K+\theta_1+\theta_2+\theta_3+1)-1)\frac{|\nabla w|^2}{w}\\
	&& +\frac{A}{e}(\ln\frac{A}{w}-1)(\theta_1+\theta_2) +\phi |\nabla q|^2 \frac{A}{e}.
	\end{eqnarray}
	For simplicity we set $\xi = 2K+\theta_1+\theta_2+\theta_3+1$. Clearly at $t=0$ we have $F(x,0)\le 0\le \frac{A}{e} (\ln\frac{A}{w}-1)(\theta_1+\theta_2)+\frac{A\theta_4^2}{e\xi},\forall x\in M$, thus in order to apply maximum principle we choose $\phi$ so that $\phi_t + \phi\xi-1\le 0$. This will give $\displaystyle F\le \frac{A}{e}\left( (\ln\frac{A}{w}-1)(\theta_1+\theta_2)+\frac{\theta_4^2}{\xi}\right)$ and hence our claim (\ref{eqTheoremLvTh}) will be settled. What remains is the existence of such function $\phi$. For this we compute
	\begin{eqnarray}
	\nonumber \phi_t + \phi\xi-1 \le 0 \implies \frac{d\phi}{1-\xi\phi} \le dt.
	\end{eqnarray}
	Integrating on $[0,t]$ we get
	\begin{eqnarray}
	\phi &\le& \frac{e^{(2K+\theta_1+\theta_2+\theta_3+1)t}-1}{(2K+\theta_1+\theta_2+\theta_3+1)e^{(2K+\theta_1+\theta_2+\theta_3+1)t}},
	\end{eqnarray}
	where we used $\xi = 2K+\theta_1+\theta_2+\theta_3+1$. Recall that $e^{\mu t}-1\ge \mu t$, consequently $\frac{t}{1+\mu t}\le \frac{e^{\mu t}-1}{\mu e^{\mu t}}$. Set $\mu=2K+\theta_1+\theta_2+\theta_3+1$ and choose $\phi = \frac{t}{1+(2K+\theta_1+\theta_2+\theta_3+1)t}$, we see all the criteria has been satisfied. This completes the proof.
\end{proof}
\begin{remark}\label{remark1}
	It should be mentioned here that the extra `$+1$' in $\xi$ [see equation (\ref{eqTheoremLvTh})] appeared due to the presence of $q(x,t)$ in our equation (\ref{eq_heat_old}). If one considers $q(x,t)$ to be a constant function then there will be no need to deduce (\ref{eq_Youngcauchy}) hence there will be no appearance of `$+1$'.
\end{remark}
\begin{corollary}\label{corr_hamilton_global}
	On a compact weighted Riemannian manifold $M$ with $Ric_f\ge - Kg$, $K\ge 0$. If $w$ is a positive solution of $\mathcal{L}_f w = qw+G(w)$, where $|q|\le\theta_1,\  |G(w)|\le\theta_2 |w|,\ |G'(w)|\le\theta_3$ then
	\begin{eqnarray}\label{eqCorrLvTT}
	 \frac{|\nabla w|^2}{w} &\le& \frac{A}{e} \left((\ln\frac{A}{w}-1)(\theta_1+\theta_2)+\ln(\frac{A}{w}) + \frac{\theta_4^2}{\xi}\right)\xi,
	\end{eqnarray}
	where $\xi = 2K + \theta_1+\theta_2 + \theta_3 + 1$. In particular if $\theta_1=\theta_2=0$, then any positive solution of $(\mathcal{L}_f-q-\partial_t)w = G(w)$, on a compact weighted Riemannian manifold $M$ with $Ric_f\ge0$ is constant.
\end{corollary}
\begin{proof}
	Letting $t\to\infty$ in (\ref{eqTheoremLvTh}) we get (\ref{eqCorrLvTT}).\\
Observe that if $M$ is a weighted Riemannian manifold with $Ric_f\ge 0$ then any positive solution $w$ of $\mathcal{L}_f w = qw+G(w)$, where $|q|\le\theta_1,\ |G(w)|\le\theta_2 |w|,\ |G'(w)|\le \theta_3$, satisfies
\begin{equation*}
 \frac{|\nabla w|^2}{w} \le \frac{A}{e} \left((\ln\frac{A}{w}-1)(\theta_1+\theta_2)+\ln(\frac{A}{w}) + \frac{\theta_4^2}{\theta_1+\theta_2 + \theta_3 + 1}\right)\left( \theta_1+\theta_2 + \theta_3 + 1\right).
\end{equation*}
The only condition required for $|\nabla w|$ to be zero is that $\theta_1=\theta_2=0$ as then $\theta_3,\theta_4$ will be zero trivially, consequently the `$+1$' term will vanish as mentioned in Remark \ref{remark1}.
\end{proof}
\begin{remark}
	The condition $\theta_1=\theta_2=0$, implies that $q\equiv 0$ and $G(w)\equiv 0$ which reduces exactly to the case of \cite[Corollary~5.8]{Yi-Li}.
\end{remark}
The above corollary is the Liouville type theorem, in which one finds certain restrictions on bounds, curvature etc. for which a solution to some PDE will be constant.
\subsection{Hessian estimate:}\label{sec3}
Before going into the main results first we define some positive constants $K_i,\ i=1,2,\cdots 8$ and some restrictions as follows
\begin{eqnarray}
	\nonumber |Rm|\le K_1,\ |\nabla Rm|\le K_2, \ |\nabla q|\le K_3, \ |\nabla^2q|\le K_4,\\
	\nonumber |\nabla^2G(w)|\le K_5 |w|, \ |q|\le K_6, \ |G(w)|\le K_7|w|, \ |\nabla G(w)|\le K_8 |w|.
\end{eqnarray}
These assumptions are necessary to derive the estimation.
Next we consider a compact domain 
$$Q_{R,T} = B(x_0,R)\times [t_0-T,t_0]\subset M\times(-\infty,\infty),$$ where $B(x_0,R)$ is a geodesic ball of radius $R$ centered at some point $x_0\in M$, $T>0$, $0\le t_0\le T$. Following \cite{Han-Zhang,Wang-hessian} we consider a function 
\begin{eqnarray}
	\label{eq_H} H = \frac{|\nabla^2 w|}{w}+\beta\frac{|\nabla w|^2}{w^2}, 
\end{eqnarray}
where $\beta$ is a positive constant to be chosen later.
\begin{lemma}\label{lemma_hess1}
	Let $(M^n,g,e^{-f}d\mu)$ be a complete weighted Riemannian manifold of dimension $n\ge 2$ with $$Ric_f\ge-K_1g,\ \ \ |Rm|\le K_1,$$ in $Q_{\frac{R}{2},\frac{T}{2}}$. If $w$ is a positive solution of (\ref{eq_heat_old}) then for any $(x,t)\in Q_{\frac{R}{2},\frac{T}{2}}$, $0<\delta<1$ and $\beta\ge \sqrt{\frac{\delta}{1-\delta}}$, the function $H$ satisfies
	\begin{eqnarray}\label{eq_lemma_hess1}
		\nonumber (\mathcal{L}_f-\partial_t)H &\ge& -2\langle\nabla H,\nabla \log w\rangle - (\Omega+\Lambda) H -\frac{\Omega}{\beta} -K_4 - K_5 + 2\beta \delta H^2 \\
		&&- 4\delta \beta^2 H\frac{|\nabla w|^2}{w^2},
	\end{eqnarray}
where $\Lambda = \max\{2K_6+CK_1+K_7, 2K_1+4K_6+4K_8\}$, $\Omega = 2K_3 + CK_2+2\beta K_3$ for some positive constant $C$.
\end{lemma}
\begin{proof}
	Let $w$ be a positive solution of (\ref{eq_heat_old}) on $M$. Using Lemma \ref{lemma_hessian_prereq} we compute
	\begin{eqnarray}
		\nonumber \mathcal{L}_f|\nabla^2 w|^2 &=& \Delta |\nabla^2 w|^2-\langle \nabla f,\nabla|\nabla^2w|^2 \rangle\\
		\nonumber &=& \sum_{i,j,k}\left[(w^2_{ij})_{kk}-f_k(w^2_{ij})_{kk}\right]\\
		\nonumber &=& \sum_{i,j,k}\left[2w_{ij}(w_{ijkk}-w_{ijk}f_k)+2w^2_{ijk}\right]\\
		\nonumber &=& 2|\nabla^2 w|\mathcal{L}_f|\nabla^2 w|+2|\nabla^3 w|^2.
	\end{eqnarray}
Since $|\nabla^3w|\ge |\nabla|\nabla^2 w||$ so the above equation reduces to
 \begin{eqnarray}
 	\nonumber \mathcal{L}_f|\nabla^2 w|^2 &\ge& 2|\nabla^2 w|\mathcal{L}_f|\nabla^2 w|+2|\nabla|\nabla^2 w||^2.
 \end{eqnarray}
By Ricci identity we find that
\begin{eqnarray}
	\nonumber w_{ijkk} &=& w_{kkij}+\sum_{l}(R_{kikl,i}w_l+R_{kijl,k}w_l+R_{kikl}w_{lj}+R_{kjkl}w_{li}+R_{kijl}w_{kl}).
\end{eqnarray}
Since $|Rm|\le K_1$ hence Shi's \cite{Shi} local derivative estimate gives the existence of  a constant $K_2$ depending on $K_1$ and $n$ such that $|\nabla Rm|\le K_2$. As in \cite{JYLI,Liu,JSUN} using our curvature restrictions we deduce
\begin{eqnarray}\label{eq_hess2.51}
	\nonumber \mathcal{L}_f |\nabla^2 w| &\ge& \frac{\langle\nabla^2 w,\nabla^2(\mathcal{L}_f w)\rangle}{|\nabla^2 w|}-CK_1|\nabla^2 w| - CK_2|\nabla w| \\
	&& + \frac{\langle \nabla^2 w,\nabla^2\langle \nabla f,\nabla w\rangle\rangle}{|\nabla^2 w|}-\frac{\displaystyle\sum_{i,j,k} w_{ij} f_k w_{ijk}}{|\nabla^2 w|}.
\end{eqnarray}
In similar way we find that
\begin{eqnarray}
	\partial_t |\nabla^2 w| &=& \frac{\langle \nabla^2 w, \nabla^2 w_t \rangle}{|\nabla^2w|}.
\end{eqnarray}
Combining these two equations
\begin{eqnarray}\label{eq_hess2.53}
	\nonumber (\mathcal{L}_f-\partial_t)|\nabla^2 w| &\ge& \frac{\langle\nabla^2 w,\nabla^2\{(\mathcal{L}_f-\partial_t) w\}\rangle}{|\nabla^2 w|}-CK_1|\nabla^2 w| - CK_2|\nabla w|.\\
\end{eqnarray}
Here the last two terms of (\ref{eq_hess2.51}) are dropped due to the following calculation.
\begin{eqnarray}
	\nonumber \nabla^2 w \langle \nabla f, \nabla w \rangle &=& \sum_{i,j,k}(f_k w_k)_{ij}\\
	\nonumber &=& \sum_{i,j,k}(f_{ki} w_k + f_k w_{ki})_{j}\\
	\nonumber &=& \sum_{i,j,k}(f_{kij} w_k + f_{ki} w_{kj} + f_{kj} w_{ki} + f_k w_{ijk})\\
	\nonumber &=& \langle \nabla^3 f,\nabla w \rangle + 2\langle \nabla^2 f,\nabla^2 w \rangle +\sum_{i,j,k} f_k w_{ijk}\\
	\nonumber &\ge& \sum_{i,j,k} f_k w_{ijk},
\end{eqnarray}
where non-negativity of the metric and commutativity of partial derivatives have been used. Using (\ref{eq_heat_old}) in (\ref{eq_hess2.53}) we get
\begin{eqnarray}
	\nonumber (\mathcal{L}_f-\partial_t)|\nabla^2 w| &\ge& \frac{\langle\nabla^2 w,\nabla^2\{qw+G(w)\}\rangle}{|\nabla^2 w|}-CK_1|\nabla^2 w| - CK_2|\nabla w|.\\
\end{eqnarray}
By Cauchy-Schwarz inequality and the bounds of $q,\ G,\ \nabla q,\ \nabla G,\ \nabla^2q,\ \nabla^2G$ we find that
\begin{eqnarray}\label{eq_hess_evol1}
	\nonumber (\mathcal{L}_f-\partial_t)|\nabla^2 w| &\ge& -(2K_3+CK_2)|\nabla w| - (K_6+CK_1)|\nabla^2 w| - K_4|w| -K_5|w|.\\
\end{eqnarray}
A direct computation gives
\begin{eqnarray}
	\nonumber (\mathcal{L}_f-\partial_t)\left(\frac{|\nabla^2 w|}{w}\right) &=& -\frac{2}{w}\left\langle \nabla (\frac{|\nabla^2 w|}{w}),\nabla w \right\rangle + \frac{1}{w} (\mathcal{L}_f-\partial_t)|\nabla^2 w|\\
	&& - \frac{|\nabla^2 w|}{w^2}(\mathcal{L}_f-\partial_t) w.
\end{eqnarray}
Using (\ref{eq_heat_old}) and (\ref{eq_hess_evol1}) in the above equation we get
\begin{eqnarray}\label{eq_hess_2.57}
	\nonumber (\mathcal{L}_f-\partial_t)\frac{|\nabla^2 w|}{w} &\ge& -\frac{2}{w}\left\langle \nabla(\frac{|\nabla^2 w|}{w}),\nabla w \right\rangle -(2K_3+CK_2)\frac{|\nabla w|}{w} \\
	&&-(2K_6+CK_1+K_7) \frac{|\nabla^2 w|}{w} -K_4 -K_5.
\end{eqnarray}
Similarly using Lemma \ref{lemma_hessian_prereq} we deduce
\begin{eqnarray}\label{eq_hess_A}
	\nonumber (\mathcal{L}_f-\partial_t)\frac{|\nabla w|^2}{w^2} &=& \frac{1}{w^2}(\mathcal{L}_f-\partial_t)|\nabla w|^2 -\frac{2}{w^3}|\nabla w|^2 (\mathcal{L}_f-\partial_t)w + \frac{6}{w^4}|\nabla w|^4\\
	&&-\frac{4}{w^3}\langle \nabla w,\nabla|\nabla w|^2 \rangle.
\end{eqnarray}
Using Lemma \ref{lemma_Bochner} (weighted Bochner formula) and (\ref{eq_heat_old}) in (\ref{eq_hess_A}) we infer
\begin{eqnarray}
	\nonumber (\mathcal{L}_f-\partial_t)\frac{|\nabla w|^2}{w^2} &=& \frac{2}{w^2} |\nabla^2 w|^2 + \frac{2}{w^2}\langle \nabla (qw+G(w)),\nabla w \rangle + \frac{2}{w^2}Ric_f(\nabla w,\nabla w)\\
	\nonumber &&- \frac{2}{w^3} |\nabla w|^2(qw+G(w)) + \frac{6}{w^4}|\nabla w|^4 - \frac{8}{w^3} \nabla^2w(\nabla w,\nabla w).\\
\end{eqnarray}
Applying bounds for the quantities and Cauchy-Schwarz inequality we get
\begin{eqnarray}\label{eq_hess2.60}
	\nonumber (\mathcal{L}_f-\partial_t)\frac{|\nabla w|^2}{w^2} &\ge& \frac{2}{w^2} \sum_{i,j} w_{ij}^2 - 2K_3\frac{|\nabla w|}{w} -(2K_6+2K_8+2K_1)\frac{|\nabla w|^2}{w^2}\\
	&& + \frac{6}{w^4} \sum_{i,j} w_i^2w_j^2 - \frac{8}{w^3}\sum_{i,j}w_{ij}w_iw_j.
\end{eqnarray}
By Young's inequality we get
\begin{eqnarray}\label{eq_hess_Young1}
	\frac{4\sum_{i,j}w_{ij}w_iw_j}{w^3} &\le& 2(1-\delta)\frac{\sum_{i,j}w_{ij}^2}{w^2} + \frac{2}{1-\delta}\left(\frac{\sum_{i}w_i^2}{w^2}\right)^2,
\end{eqnarray}
for some $0<\delta<1$. Also it can be easily seen that
\begin{eqnarray}\label{eq_hess_Young2}
	\nonumber && -\frac{8}{w^3}\sum_{i,j}w_{ij}w_i w_j + \frac{6}{w^4} \sum_{i,j}w_i^2w_j^2\\
	&& = -2\left\langle \nabla(\frac{|\nabla w|^2}{w^2}),\nabla\log w \right\rangle  - \frac{4}{w^3}\nabla^2w(\nabla w,\nabla w) + \frac{2}{w^4}|\nabla w|^4.
\end{eqnarray}
Using Young's inequality, (\ref{eq_hess_Young1}), (\ref{eq_hess_Young2}) and a scaling with constant $\beta>0$ in (\ref{eq_hess2.60}) we obtain
\begin{eqnarray}
	\nonumber (\mathcal{L}_f-\partial_t) \left(\beta\frac{|\nabla w|^2}{w^2}\right) &\ge& 2\beta\delta \frac{\sum_{i,j}w_{ij}^2}{w^2} - 2\beta^3\frac{|\nabla w|}{w} - \beta (2K_1+4K_6+4K_8)\\
	\nonumber && -2\left\langle \nabla(\beta\frac{|\nabla w|^2}{w^2}),\nabla \log w \right\rangle - \frac{2\beta}{1-\delta}\left(\frac{\sum_{i}w_i^2}{w^2}\right)^2 \\
	&&+ \frac{2\beta}{w^4}\sum_{i,j}w_i^2w_j^2.
\end{eqnarray}
Combining the above equation with (\ref{eq_hess_2.57}) we find that the function $H$ satisfies
\begin{eqnarray}
	\nonumber (\mathcal{L}_f-\partial_t)H &\ge& -2\langle \nabla H,\nabla\log w \rangle - (2K_3+CK_2+2\beta K_3)\frac{|\nabla w|}{w} -\Lambda H - K_4 - K_5\\
	&&+2\beta\delta \frac{|\nabla^2w|^2}{w^2} - \frac{2\beta}{1-\delta}\frac{|\nabla w|^4}{w^4} + 2\beta \frac{|\nabla w|^4}{w^4}.
\end{eqnarray}
Given that $\frac{|\nabla^2 w|}{w} = H - \beta\frac{|\nabla w|^2}{w^2}$ and $H\ge \beta \frac{|\nabla w|^2}{w^2}$ so we infer
\begin{eqnarray}
	\nonumber (\mathcal{L}_f-\partial_t)H &\ge& -2\langle \nabla H,\nabla\log w \rangle - (2K_3+CK_2+2\beta K_3)\frac{|\nabla w|}{w} -\Lambda H - K_4 - K_5\\
	&&+2\beta\delta H^2 - 4H\beta^2\delta\frac{|\nabla w|^2}{w^2} + (2\beta^3+2\beta-\frac{2\beta}{1-\delta}) \frac{|\nabla w|^4}{w^4}.
\end{eqnarray}
Since $\beta\ge\sqrt{\frac{\delta}{1-\delta}}>0$ hence
	\begin{eqnarray}
		\nonumber (\mathcal{L}_f-\partial_t)H &\ge& -2\langle \nabla H,\nabla\log w \rangle - \Omega\frac{|\nabla w|}{w} -\Lambda H - K_4 - K_5 \\
		&&+2\beta\delta H^2 - 4H\beta^2\delta\frac{|\nabla w|^2}{w^2}.
	\end{eqnarray}
Applying an elementary inequality $-\Omega \frac{|\nabla w|}{w}\ge -\Omega H -\frac{\Omega}{\beta}$ in the above equation gives (\ref{eq_lemma_hess1}). This completes the proof.
\end{proof}
In the next result we are going to derive bound for $\frac{|\nabla^2 u|}{u}$ using the above lemma by considering a cut-off function. An important note regarding the above theorem is made at the end of the proof.
\begin{theorem}\label{theorem_hess1}
	Let $(M^n,g,e^{-f}d\mu)$ be a complete weighted Riemannian manifold of dimension $n\ge 2$ with $$Ric_f\ge-K_1g,\ \ \ |Rm|\le K_1,$$ in $Q_{\frac{R}{2},\frac{T}{2}}$. If $w$ is a positive solution of (\ref{eq_heat_old}) then for any $(x,t)\in Q_{\frac{R}{2},\frac{T}{2}}$, $0<\delta<1$ and $\beta\ge \sqrt{\frac{\delta}{1-\delta}}$, we have
	\begin{eqnarray}
		\nonumber \frac{|\nabla^2 w|}{w} &\le& (4\sqrt2 -1)\beta\frac{|\nabla w|^2}{w^2}+\sqrt2 \Big(\frac{C}{\delta^4\beta^6R^4}+\frac{\Omega}{\delta\beta^2}+\frac{K_4+K_5}{\delta\beta}+\frac{C}{\delta^2 \beta^2 T^2}\\
		&&+\frac{C}{\delta^2\beta^2}(\frac{1}{R^4}+\frac{K_1}{R^2})\Big)^\frac12+\sqrt2 \Big(\frac{\Omega+\Lambda}{\delta\beta}+\frac{2C^2}{\delta\beta T^2}\Big),
	\end{eqnarray}
where $\Lambda = \max\{2K_6+CK_1+K_7, 2K_1+4K_6+4K_8\}$, $\Omega = 2K_3 + CK_2+2\beta K_3$ for some positive constant $C$.
\end{theorem}
\begin{proof}
	With the help of \cite{Bali-Cao-Pule,SOUPLET-ZHANG} and using geodesic polar coordinates we consider a cutoff function $\psi\equiv\psi(x,t)$ supported in $Q_{R,T}$ satisfying
	\begin{enumerate}
		\item $\psi(r,t)=\psi(d(x,x_0),t),\ \psi(r,t)=1$ in $Q_{\frac{R}{2},\frac{T}{2}},\ 0\le\psi\le 1$.
		\item $\frac{\partial_r\psi}{\psi^\alpha}\le\frac{C_\alpha}{R},\ \frac{|\partial^2_r \psi|}{\psi^\alpha}\le \frac{C_\alpha}{R^2}$ with $0<\alpha<1$.
		\item  $\frac{\partial_t \psi}{\sqrt\psi}\le\frac{C}{T}$.
		\item  $\psi$ is radially decreasing in the spatial variables,
	\end{enumerate}
$C_\alpha, C$ are positive constants. A direct calculation using Lemma \ref{lemma_hessian_prereq} gives
\begin{eqnarray}
	(\mathcal{L}_f-\partial_t)(\psi H) &=& \psi(\mathcal{L}_f-\partial_t)H + H(\mathcal{L}_f-\partial_t)\psi + 2\langle \nabla\psi,\nabla H \rangle.
\end{eqnarray}
Applying Lemma \ref{lemma_hess1} we infer
\begin{eqnarray}
	\nonumber (\mathcal{L}_f-\partial_t)(\psi H) &=& \psi\Big[-2\langle \nabla H,\nabla\log w \rangle - (\Omega+\Lambda)H -\frac{\Omega}{\beta} - K_4 - K_5 + 2\beta\delta H^2 \\
	&& - 4\delta\beta^2 H\frac{|\nabla w|^2}{w^2}
	\Big] + H(\mathcal{L}_f-\partial_t)\psi + 2\langle \nabla\psi,\nabla H \rangle
\end{eqnarray}
Suppose that the function $\psi H$ achieves its maximum at a point $(x_1,t_1)$, without loss of generality which is not on the cut-locus of $M$ \cite{CALABI}. Hence at $(x_1,t_1)$ we find
\begin{eqnarray}
	\begin{cases}
		\nonumber (\mathcal{L}_f-\partial_t)(\psi H) \le 0,\\
		\nonumber \nabla(\psi H) = 0, \text{i.e., } \psi \nabla H=-H\nabla\psi.
	\end{cases}
\end{eqnarray}
Thus we have
\begin{eqnarray}
	\nonumber0 &\ge& \psi\left[-2\langle \nabla H, \nabla \log w \rangle - (\Omega+\Lambda)H -\frac{\Omega}{\beta} - K_4 - K_5 +2\beta\delta H^2 - 4\delta\beta^2 H\frac{|\nabla w|^2}{w^2}\right]\\
	&&+H(\mathcal{L}_f-\partial_t)\psi-2H\frac{|\nabla\psi|^2}{\psi}.
\end{eqnarray}
Dividing the above equation by $\delta\beta$ and using the following results
\begin{eqnarray}
	\frac{2}{\delta\beta}H\langle\nabla\psi,\nabla \log w\rangle &\le& \frac{1}{3} \psi H^2 + \frac{c}{\delta^4\beta^6R^4},\\
	\frac{1}{\delta\beta}H\Delta\psi &\le& \frac{1}{3} \psi H^2 +\frac{C}{\delta^2\beta^2}\left(\frac{1}{R^4}+\frac{K_1}{R^2}\right),\\
	\frac{1}{\delta\beta}H\partial_t \psi &\le& \frac{1}{3}\psi H^2 +\frac{C}{\delta^2\beta^2T^2},
\end{eqnarray}
we get
\begin{eqnarray}\label{eq_hess_2.74}
	\nonumber \psi H^2 &\le& \Big( \frac{C}{\delta^4 \beta^6 R^4}+\frac{\Omega \psi}{\delta\beta^2}+\frac{(K_4+K_5)\psi}{\delta\beta}+\frac{C}{\delta^2\beta^2 T^2} + \frac{C}{\delta^2\beta^2}(\frac{1}{R^4}+\frac{K_1}{R^2}) \Big)\\
	&&+\Big(\frac{\psi (\Omega+\Lambda)}{\delta\beta}+4\psi\beta\frac{|\nabla w|^2}{w^2}+\frac{2C^2}{\delta\beta T^2}\Big)H.
\end{eqnarray}
Using Young's inequality we have
\begin{eqnarray}\label{eq_hess_2.75}
	\nonumber (A+4\psi\beta\frac{|\nabla w|^2}{w^2})H &\le& \frac{A^2}{\psi} + \frac{\psi H^2}{4} + 16\beta^2 \frac{|\nabla w|^4}{w^4}+\frac{\psi H^2}{4}\\
	&\le& \psi\left( \frac{\Omega+\Lambda}{\delta\beta}+\frac{2C^2}{\delta\beta T^2\psi} \right)^2 + 16\beta^2\frac{|\nabla w|^4}{w^4} + \frac{\psi H^2}{2},
\end{eqnarray}
where $A=\frac{\psi(\Omega+\Lambda)}{\delta\beta}+\frac{2C^2}{\delta\beta T^2} $. Putting (\ref{eq_hess_2.75}) in (\ref{eq_hess_2.74}) and using the inequality we get
\begin{eqnarray}
	\psi H^2 &\le& 2\mathcal{A}+2\psi \mathcal{B}^2+32\beta^2\frac{|\nabla w|^4}{w^4},
\end{eqnarray}
where $$\mathcal{A}=\frac{C}{\delta^4\beta^6 R^4}+\frac{\Omega\psi}{\delta\beta^2}+\frac{(K_4+K_5)\psi}{\delta\beta}+\frac{C}{\delta^2\beta^2T^2}+\frac{C}{\delta^2\beta^2}(\frac{1}{R^4}+\frac{K_1}{R^2}),$$
$$\mathcal{B} = \frac{\Omega+\Lambda}{\delta\beta}+\frac{2C^2}{\delta\beta T^2\psi}.$$
For all $(x,t)\in Q_{R,T}$ we have
$$(\psi H)^2(x,t)\le(\psi H)^2(x_1,t_1)\le\psi H^2(x_1,t_1).$$
Also recall that $\psi\equiv 1$ in $Q_{\frac{R}{2},\frac{T}{2}}$. Thus by definition of $H$ we deduce 
\begin{eqnarray}
	\nonumber \frac{|\nabla^2 w|}{w} &\le& (4\sqrt2 -1)\beta\frac{|\nabla w|^2}{w^2}+\sqrt2 \Big(\frac{C}{\delta^4\beta^6R^4}+\frac{\Omega}{\delta\beta^2}+\frac{K_4+K_5}{\delta\beta}+\frac{C}{\delta^2 \beta^2 T^2}\\
	&&+\frac{C}{\delta^2\beta^2}(\frac{1}{R^4}+\frac{K_1}{R^2})\Big)^\frac12+\sqrt2 \Big(\frac{\Omega+\Lambda}{\delta\beta}+\frac{2C^2}{\delta\beta T^2}\Big),
\end{eqnarray}
where the inequality $\sqrt{x+y}\le\sqrt{x}+\sqrt{y}$ have been used. This completes the proof.
\end{proof}
\begin{remark}\label{remark_hess1}
	Note that in the above theorem we have not put any restriction on the quantity $\frac{|\nabla w|}{w}$ which represents the gradient estimate for the solution $w$ of (\ref{eq_heat_old}). The idea behind this kind of extension is to get various Hessian estimations from different kinds of gradient estimates. 
\end{remark}
As an immediate consequence we have the global Hessian estimate as follows.
\begin{corollary}\label{corr_hess2.4}
	Let $(M^n,g,e^{-f}d\mu)$ be a complete weighted Riemannian manifold of dimension $n\ge 2$ with $$Ric_f\ge-K_1g,\ \ \ |Rm|\le K_1,$$ on $M\times [t_0-T,t_0]$. If $w$ is a positive solution of (\ref{eq_heat_old}) then for any $(x,t)\in M\times[t_0-T,t_0]$, $0<\delta<1$ and $\beta\ge \sqrt{\frac{\delta}{1-\delta}}$, we have
	\begin{eqnarray}
		\nonumber \frac{|\nabla^2 w|}{w} &\le& (4\sqrt2 -1)\beta\frac{|\nabla w|^2}{w^2}+\sqrt2 \Big(\frac{\Omega}{\delta\beta^2}+\frac{K_4+K_5}{\delta\beta}+\frac{C}{\delta^2 \beta^2 T^2}\Big)^\frac12\\
		&&+\sqrt2 \Big(\frac{\Omega+\Lambda}{\delta\beta}+\frac{2C^2}{\delta\beta T^2}\Big),
	\end{eqnarray}
	where $\Lambda = \max\{2K_6+CK_1+K_7, 2K_1+4K_6+4K_8\}$, $\Omega = 2K_3 + CK_2+2\beta K_3$ for some positive constant $C$.
\end{corollary}
\begin{proof}
	Letting $R\to+\infty$ completes the proof.
\end{proof}
Recalling Remark \ref{remark_hess1} we now show how one can derive different Hessian estimates using known gradient estimates. These corollaries are basically applications of Theorem \ref{theorem_hess1}.
\begin{corollary}\label{corr_LYHess}
In view of (\ref{eq_lambda}) with the assumptions of Corollary \ref{corr1} and Corollary \ref{corr_hess2.4} we have the following Li-Yau type Hessian estimate on $M\times (0,T]$ given by
\begin{eqnarray}\label{eq_corr_LYHess}
	\nonumber\frac{|\nabla^2w|}{w} &\le& (4\sqrt2 -1)\beta \left(\Lambda_{\alpha,\epsilon}+\frac{m\alpha^2}{2t(1-\epsilon)}+\alpha(\frac{w_t}{w}+q(x,t)+\frac{G(w)}{w})\right)\\
	&& +\sqrt2 \mathcal{A}^\frac12 + \sqrt2 \mathcal{B},
\end{eqnarray}
where $\displaystyle \mathcal{A}=\frac{\Omega}{\delta\beta^2}+\frac{K_4+K_5}{\delta\beta}+\frac{C}{\delta^2 \beta^2 T^2}$ and $\displaystyle \mathcal{B}=\frac{\Omega+\Lambda}{\delta\beta}+\frac{2C^2}{\delta\beta T^2}$, rest of the symbols are defined earlier.
\end{corollary}
\begin{proof}
	Set $u=\log w$ in Corollary \ref{corr1} and using (\ref{eq_lambda}) (to cover all the cases as defined in Corollary \ref{corr1}) we find an inequality related to $\frac{|\nabla w|^2}{w^2}$, which is then used in Corollary \ref{corr_hess2.4} to get (\ref{eq_corr_LYHess})
\end{proof}
Using the above Corollary we derive a local time reverse Harnack inequality.
\begin{corollary}\label{corr_revHarnack}
	For any $x\in Q_{R,T}$ and $t_1,t_2\in (0,T]$ with $t_1<t_2$ if $w$ is a solution of (\ref{eq_heat_old}) satisfying the conditions of Corollary \ref{corr_LYHess} then we have
	\begin{eqnarray}\label{eq_LYHess_revHarnack}
		w(x,t_2) &\le& w(x,t_1) \left(\frac{t_2}{t_1}\right)^{N_1} \exp\left\{N_2 (t_2-t_1)\right\},
	\end{eqnarray}
where 
\begin{eqnarray}
	\nonumber N_1 &=& \frac{1}{1-(4\sqrt2-1)\beta\alpha}\left\{ K_6+K_7+(4\sqrt2-1)\beta\Lambda_{\alpha,\epsilon}+\sqrt2 \mathcal{A}^\frac12 +\sqrt2 \mathcal{B}\right\},\\
	\nonumber N_2 &=& \frac{(4\sqrt2-1)\beta m\alpha^2}{2(1-\epsilon)(1-(4\sqrt2 -1)\beta\alpha)},
\end{eqnarray}
with the assumption $$0<\delta\le \frac{1}{1+(4\sqrt2-1)^2\alpha^2}\text{ so that }\sqrt{\frac{\delta}{1-\delta}}\le\beta\le \frac{1}{(4\sqrt2-1)\alpha},$$
where $$\mathcal{A}=\frac{C}{\delta^4\beta^6 R^4}+\frac{\Omega}{\delta\beta^2}+\frac{(K_4+K_5)}{\delta\beta}+\frac{C}{\delta^2\beta^2T^2}+\frac{C}{\delta^2\beta^2}(\frac{1}{R^4}+\frac{K_1}{R^2}),$$
$$\mathcal{B} = \frac{\Omega+\Lambda}{\delta\beta}+\frac{2C^2}{\delta\beta T^2},$$ and rest of the symbols are defined earlier.
\end{corollary}
\begin{proof}
	Rewriting (\ref{eq_corr_LYHess}) in local form we find
	\begin{eqnarray}
		\nonumber\frac{w_{ij}}{w} &\le& (4\sqrt2 -1)\beta \left(\Lambda_{\alpha,\epsilon}+\frac{m\alpha^2}{2t(1-\epsilon)}+\alpha(\frac{w_t}{w}+q(x,t)+\frac{G(w)}{w})\right)\\
		&& +\sqrt2 \mathcal{A}^\frac12 + \sqrt2 \mathcal{B}.
	\end{eqnarray}
Following \cite[Theorem~5.1]{Wang-hessian} and using the non-negativity of the metric we infer
\begin{eqnarray}
	\nonumber \frac{\mathcal{L}_f w}{w} &\le& (4\sqrt2 -1)\beta \left(\Lambda_{\alpha,\epsilon}+\frac{m\alpha^2}{2t(1-\epsilon)}+\alpha(\frac{w_t}{w}+q(x,t)+\frac{G(w)}{w})\right)\\
	&& +\sqrt2 \mathcal{A}^\frac12 + \sqrt2 \mathcal{B}.
\end{eqnarray}
Applying (\ref{eq_heat_old}) and invoking the restrictions on $\beta, \delta$ we deduce
\begin{eqnarray}
	\partial_t(\log w) &\le& N_1 + \frac{N_2}{t}.
\end{eqnarray}
Integrating between the time $t_1$ and $t_2$ we get
\begin{eqnarray}
	\log\left(\frac{w(x,t_2)}{w(x,t_1)}\right) &\le& N_1 (t_2-t_1) + N_2 \log\left(\frac{t_2}{t_1}\right). 
\end{eqnarray}
Upon exponentiation we get (\ref{eq_LYHess_revHarnack}). This completes the proof.
\end{proof}
In similar manner we can derive a Hamilton type Hessian estimate as follows.
\begin{corollary}
	For $x\in M$, $t\in (0,T]$, $t_0=T$, $0<w\le\frac{A}{e}$, $Ric_f\ge-K_1g$ and with the assumptions of Corollary \ref{corr_hamilton_global}, Corollary \ref{corr_hess2.4} we have the following Hamilton type Hessian estimate
	\begin{eqnarray}
		\nonumber\frac{|\nabla^2w|}{w} &\le& (4\sqrt2-1)\frac{\beta A}{ew}\left\{ \log(\frac{A}{w}) + (\ln\frac{A}{w}-1)(\theta_1+\theta_2) +\frac{K_3^2}{\eta} \right\} \left(\frac{1}{t}+\eta\right)\\
		&&+\sqrt2\mathcal{A}^\frac12 +\sqrt2 \mathcal{B},
	\end{eqnarray}
where $\eta = 2K_1+K_6+K_7+\theta_3+1$, $$\mathcal{A}=\frac{\Omega}{\delta\beta^2}+\frac{(K_4+K_5)}{\delta\beta}+\frac{C}{\delta^2\beta^2T^2},\ \mathcal{B} = \frac{\Omega+\Lambda}{\delta\beta}+\frac{2C^2}{\delta\beta T^2},\ |G'(w)|\le\theta_3,$$ and the rest of the symbols are defined earlier.
\end{corollary}
\begin{proof}
	The proof is similar to Corollary \ref{corr_LYHess} and Corollary \ref{corr_revHarnack}, just use the value of $\frac{|\nabla w|^2}{w}$ from Corollary \ref{corr_hamilton_global}.
\end{proof}
\section{Concluding remark}
The methods of gradient and Hessian estimation are indeed a very powerful tool to understand the nature of the solution of a heat type equation (\ref{eq_heat_old}). As shown in the introduction, this is a very rich area of research. In this article, we derived Li-Yau type estimation, Hamilton type estimation for a positive solution of (\ref{eq_heat_old}) on a static weighted Riemannian manifold (static means the metric is not evolving under any geometric flow). Next we derived Harnack type inequality with the help of Li-Yau type estimation which gives us the information about the change of heat between two points in the time involved space $M\times [0,T]$ along a geodesic. Using the Hamilton type estimation we derive some restrictions under which any solution of (\ref{eq_heat_old}) will be constant i.e., a Liouville type theorem is stated. Finally we derived Hessian estimate for the aforementioned equation and concluded several corollaries as application. A local time reversed Harnack inequality is also derived at the end. 

\vspace{.25in}
\noindent{Shyamal Kumar Hui \\ Department of Mathematics, The University of Burdwan, Golapbag, Burdwan 713104, West Bengal, India}\\
Email: skhui@math.buruniv.ac.in\\

\noindent{Abimbola Abolarinwa\\ Department of Mathematics, University of Lagos, Akoka, Lagos State, Nigeria}\\
Email: A.Abolarinwa1@gmail.com,  aabolarinwa@unilag.edu.ng\\

\noindent{Sujit Bhattacharyya \\ Department of Mathematics, The University of Burdwan, Golapbag, Burdwan 713104, West Bengal, India}\\
Email: sujitbhattacharyya.1996@gmail.com
\end{document}